\documentclass[11pt]{article}
\usepackage{graphicx}
\usepackage{amsmath}
\usepackage{amssymb}
\usepackage{amsfonts}
\usepackage{amsthm}
\usepackage{xcolor}
\usepackage{enumitem}
\usepackage{url}
\usepackage{bm}
\numberwithin{equation}{section}
\def\TT{{\bf T}}
\def\D{{\bf D}}
\def\R{\mathbb{R}}

\def\Z{\mathbb{Z}}
\def\dive{\text{div}}

\def\x {{\bf x}}
\def\y {{\bf y}}
\def\u {{\bf u}}
\def\v {{\bf v}}
\def\w {{\bf w}}
\def\f {{\bf f}}
\def\g {{\bf g}}
\def\h {{\bf h}}
\def\n {{\bf n}}
\def\0 {{\bf 0}}
\def\i {{\bf i}}
\def\j {{\bf j}}
\def\k {{\bf k}}

\def\e {\varepsilon}

\def\endproof{\hfill $\Box$\par\vskip3mm}
\def\neweq#1{\begin{equation}\label{#1}}
\def\endeq{\end{equation}}

\def\eps{\varepsilon}
\def\ds{\displaystyle}
\newtheorem{theorem}{Theorem}[section]

\newtheorem{lemma}[theorem]{Lemma}
\newtheorem{proposition}[theorem]{Proposition}

\newtheorem{remark}[theorem]{Remark}

\title{Elasticity solution for a 3D hollow cylinder\\ axially loaded at the end faces}
\author{Mario DE MIRANDA$^\ddag$, Marta DE MIRANDA$^\ddag$, Alessio FALOCCHI$^\dag$, \\Alberto FERRERO$^*$, Luca MARININI$^\ddag$, \vspace{5mm}\\
 \small $^\ddag$ Studio De Miranda Associati, Via C. Pisacane 26, Milano, Italy\\
{\small $^\dag$ Dipartimento di Matematica - Politecnico di Milano, Milano,  Italy}\\
{\small $^*$ Dipartimento di Scienze e Innovazione Tecnologica, Università del Piemonte Orientale, Alessandria, Italy}}
\date{}
\begin{document}
\maketitle
\begin{abstract}
	Starting from an applicative problem related to the modeling of an element of a cable-stayed bridge, we compute the elasticity solution for a hollow cylinder loaded at the end faces with axial loads. We prove results of symmetry for the solution and we expand it in proper Fourier series; computing the Fourier coefficients in adapted power series, we provide the explicit solution. We consider an engineering case of study, applying the corresponding approximate formula and giving some estimates on the error committed with respect to the truncation of the series.
\end{abstract}

\section{Introduction}
In the recent years the interest of the mathematicians for engineering applications has grown more and more; this is due to an evolution of the mathematics, thanks for instance to the development of new techniques to deal with nonlinear problems and the support of automatic calculators to obtain previsions unthinkable in the past.

\begin{figure}[h]
	\centering
	\includegraphics[width=15cm]{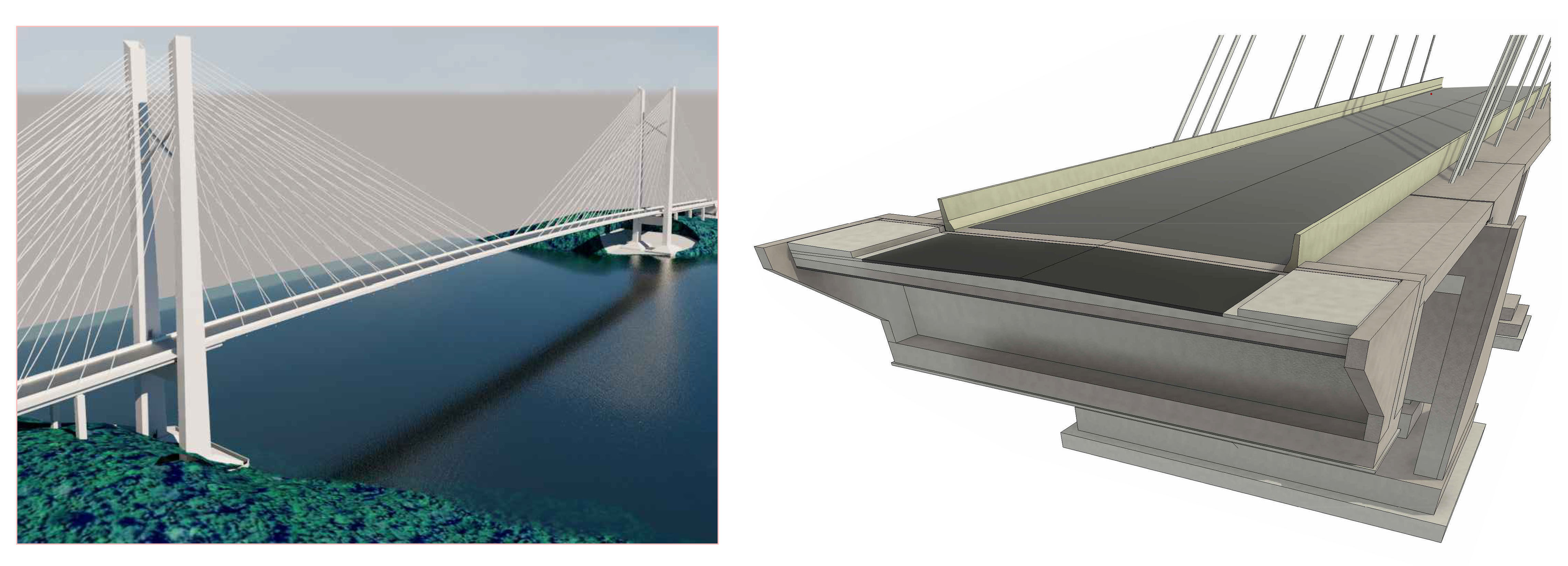}
	\caption{From the top on the left a render of a recent cable stayed bridge designed by Studio De Miranda Associati and a detail of its deck.}\label{fig:blister}
\end{figure}

The mathematical modeling of specific phenomena is one of the ambitious aims of the applied mathematics. Recently some mathematical models for suspension bridges \cite{Gazz} have been developed with the scope to understand and, then, to prevent instability phenomena; they were studied models for suspension bridges with geometrical nonlinearities, e.g. see \cite{crfaga,fa1}, models for partially hinged plates, see \cite{fega}, models for non homogeneous partially hinged plates, see \cite{befa0, befa, befafega}, models for homogeneous beams with intermediate piers \cite{gaga} and non homogeneous beams \cite{befaga}. In all these cases the application of analytical methods to real problems allowed to find suggestions and practical remedies that can be discussed with engineers.

\begin{figure}[h]
	\centering
	\includegraphics[width=15cm]{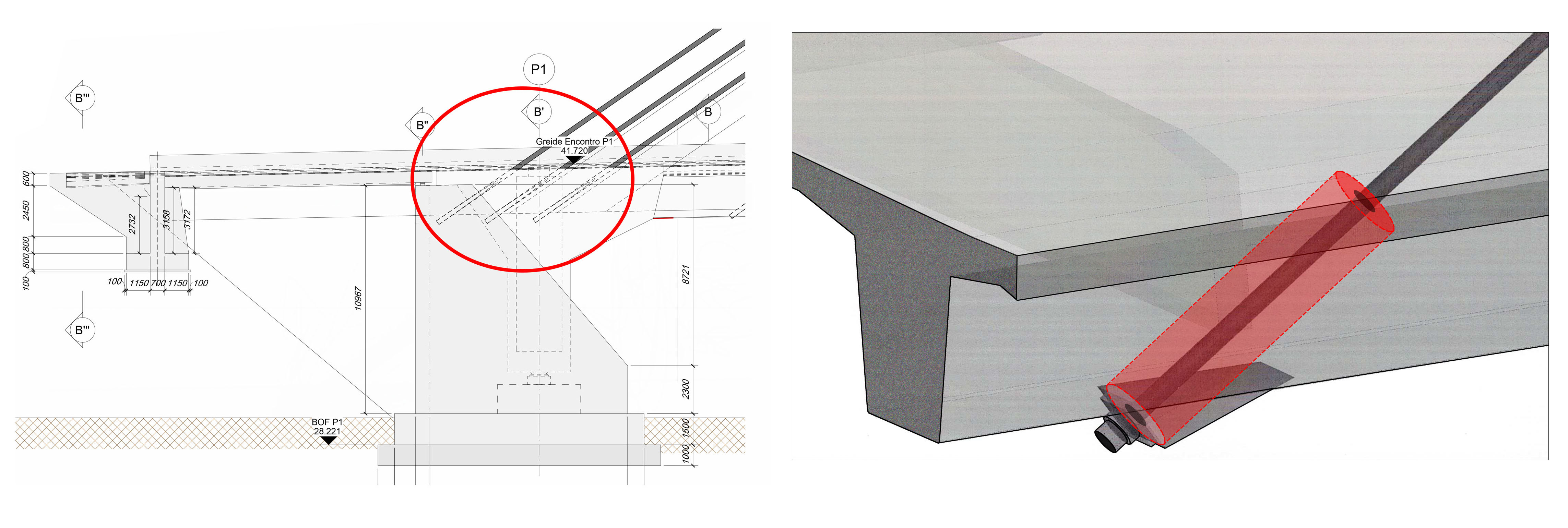}
	\caption{From the left a frontal view of blister elements and the modelization of the element through the hollow cylinder (in red).}\label{fig:blister2}
\end{figure}

 This is the aim also of this paper. Here the problem, suggested by the structural civil engineering Studio De Miranda Associati, is related to the modeling of the stresses in a constructive detail of a bridge: the blister. In the cable-stayed bridge the blister is the structural element where the steel forestay anchors to the deck. In Figure \ref{fig:blister} is shown a render of a future cable-stayed bridge, designed by Studio De Miranda Associati, that will be built in Brazil; a detail of the related blister element is given in Figure \ref{fig:blister2}.

\begin{figure}[h]
	\centering
	\includegraphics[width=17cm]{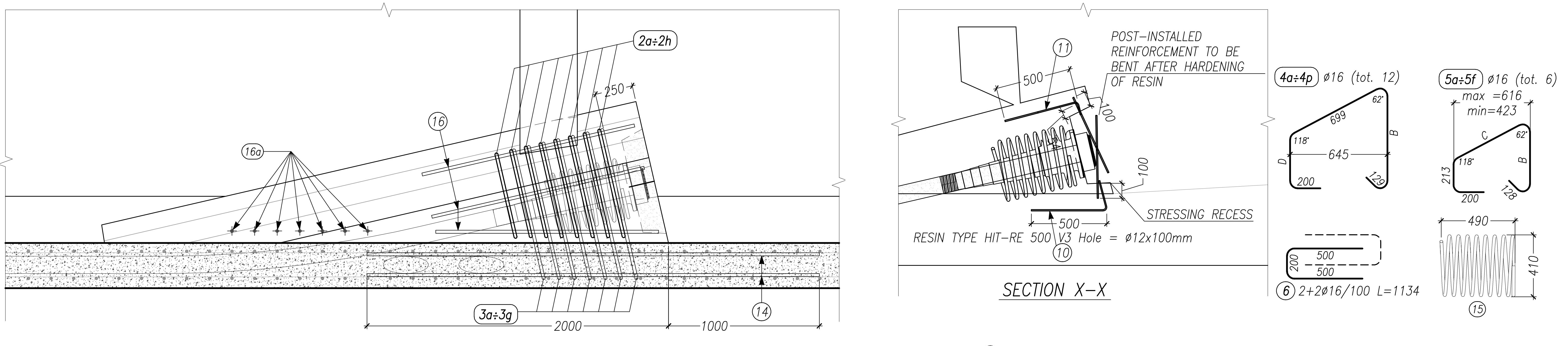}
	\caption{Detail of the executive draw of a blister of a recent prestressed concrete bridge.}\label{fig:armature}
\end{figure}

When the deck is built in reinforced concrete, as in this case, the blister is an important point to design; indeed, the high density of the steel reinforcement may cause zone with low concrete capacity and possible remarkable cracking. To have an idea of the complexity of this element a detail of the executive draw of a blister for another stayed bridge, designed by Studio De Miranda Associati, is shown in Figure \ref{fig:armature}.

For all these reasons it is important to estimate with precision the stresses acting on the element, so that the reinforcing steel in the concrete can be computed without surplus. In engineering literature some of the best known references related to the distribution of the stresses in prisms of concrete are \cite{ciria,leo}; here the authors consider many combinations of load on the prism and for each one the possible strategies to design the steel bars.
These results are obtained from particular solutions of the well known equation of the linear elasticity, see e.g. \cite{ciarlet}; we recall it here briefly in the general 3D case.

Given $\Omega\subset \R^3$ an elastic homogeneous solid body, we denote by $\u:\Omega\rightarrow\R^3$ the displacement vector at any point of the reference configuration of the elastic body itself, see the list of notations at the end of the paper.
We denote by $\TT \u$ the stress tensor and by $\lambda$ and $\mu$ the classical Lamé constants;
it is known that $\lambda$ and $\mu$ may be expressed in terms of the Young modulus $E$ and Poisson ratio $\nu\in(-1,\frac{1}{2})$ as
\begin{equation} \label{eq:lame}
  \lambda=\frac{E\nu}{(1+\nu)(1-2\nu)} \, ,  \qquad  \mu=\frac{E}{2(1+\nu)} \, .
\end{equation}
The equation of linear elasticity reads
\begin{equation}\label{pb}
\begin{cases}
-\mu \Delta \u-(\lambda+\mu)\nabla(\dive \u)=\f\quad&\text{in }\Omega, \\[7pt]
(\TT \u) \n=\g\quad &\text{on }\partial \Omega,
\end{cases}
\end{equation}
where $\f$ and $\g$ are respectively the forces per unit volume and the boundary forces per unit surface acting on $\Omega$, while $\n$ is the unit outward normal vector to $\partial\Omega$.

In Section \ref{math model} we briefly derive \eqref{pb} from variational principles and we recall the existence and uniqueness results in Theorem \ref{p:existence}; these are classical topics in linear elasticity, see e.g. \cite{ciarlet,weest}, but we recall them in our framework for completeness since the question about uniqueness of solutions of \eqref{pb} is not trivial at all and it needs an additional condition to be achieved.

The theoretical solution from which come the applicative cases considered in \cite{ciria,leo} is given in \cite{sundara}, where $\Omega$ is a rectangular prism under end loads. Thanks to this simple geometry and loading condition the authors find explicitly the solution in form of double Fourier series. The result is obtained applying the Galerkin vector method, a technique allowing to pass from the second order differential equation \eqref{pb} to a simpler biharmonic equation, see \cite{sundara}.

We point out that to find the explicit solution of \eqref{pb} for generic $\Omega$ and loading conditions is a very hard task. In this paper we find it for $\Omega$ coincident with a hollow cylinder loaded on the opposite faces, since this geometry fits the modeling of the concrete of the blister, see Figure \ref{fig:blister2} on the right; indeed, the forestay of the bridge is circular and passes through the cylindrical hole, applying a distributed load on the opposite faces due to its tensioning,  see Figure \ref{cylinder0}. 

The precise definition of the model is given in Section \ref{cylinder0}; the application of axial loads leads to a solution having axial symmetric properties, see Proposition \ref{prop0}. In the real blister it is also possible to have non radial loadings coming from the deck, but this is a first attempt of modeling that may be implemented in future works; anyway, the solution found here may have general interest beyond this specific application.

The definition of the solution is given by steps: in subsection \ref{periodic} we provide a periodic extension of the loads in the variable $z$ corresponding to the symmetry axis of the hollow cylinder, in such a way that it becomes possible to expand the solution in Fourier series with respect to the variable $z$; then we compute the Fourier coefficients which come to be functions in the other two variables $x$ and $y$, corresponding to directions orthogonal to the symmetry axis of the hollow cylinder; in subsection \ref{cil} we pass to the cylindrical coordinates and, exploiting the axial symmetry, we reduce ourself to study a system of ODEs in the radial polar coordinate $\rho$; we compute the Fourier coefficients as functions of the variable $\rho$ through an adapted expansion in power series so that we are able to state Theorem \ref{teo_soluzioni}, collecting the explicit solution.

In Section \ref{num} we give some hints to truncate the series and we apply the results to an engineering case of study. As it will be explained in details, it will be necessary to compute numerically the first $M$ terms in the Fourier series expansion with $M$ to be chosen sufficiently large in order to minimize the truncation error.
The main question in this procedure is that the computation of those Fourier coefficients, which are solutions of suitable boundary value problems of ODEs, requires the numerical resolution of some algebraic linear systems in four variables which exhibit a condition number higher and higher as $M$ grows; if we need a truncation error smaller than ours, we may consider alternative numerical procedures. We emphasize that the main purpose of this article is to obtain an analytical representation of the unique symmetric solution of \eqref{pb} in the case of the hollow cylinder with the perspective of reproducing such method in more general situations with not necessarily symmetric external loads.

As already explained in details, the main analytical and numerical results of the article are stated in Sections \ref{math model}-\ref{num} and their proofs are given in Section \ref{proofs}.
The final part of the paper is devoted to the conclusions, see Section \ref{conclusions}, and a list of notations which can be helpful for the reader.

\bigskip

\section{The definition of the mathematical model for the linear elasticity}\label{math model}
In this Section we derive the differential equations for the linear elasticity from variational arguments and we state a theorem related to the existence of solutions. Although these results are well known overall in the engineering field, we review them from a mathematical point of view, applying the Fredholm alternative to prove existence of solutions.
\subsection{The derivation of the differential equations}
We recall that $\Omega\subset \R^3$ is the domain of the elastic body and $\u$ is the displacement function with components $\u=(u_1,u_2,u_3)$. We denote by $\D \u$ the linearized strain tensor, which in the sequel will be simply called strain tensor, since we only deal with the linear theory; the stress tensor can be written as
\begin{equation} \label{eq:stress}
\TT \u=\begin{bmatrix}
\sigma^{1} & \tau^{12} & \tau^{13} \\
\tau^{12} & \sigma^{2} & \tau^{23}\\
\tau^{13} & \tau^{23} & \sigma^{3} \\
\end{bmatrix} \, .
\end{equation}
It is well known that by the Hooke's Law for isotropic materials it holds
\begin{equation}\label{Hooke}
 \TT \u=\lambda \text{tr}(\D \u)\, I+2\mu \D \u \, ,
\end{equation}
where $\lambda$ and $\mu$ are the Lamé constants.

Combining \eqref{eq:lame}, \eqref{eq:stress} and \eqref{Hooke} we infer
\begin{equation}\label{stresses}
	\begin{split}
	&\sigma^1=\dfrac{E}{(1+\nu)(1-2\nu)}\bigg[(1-\nu)\dfrac{\partial u_1}{\partial x}+\nu\bigg(\dfrac{\partial u_2}{\partial y}+\dfrac{\partial u_3}{\partial z}\bigg) \bigg]\qquad\quad \tau^{12}=\dfrac{E}{2(1+\nu)}\bigg(\dfrac{\partial u_1}{\partial y}+\dfrac{\partial u_2}{\partial x}\bigg)\\
	&\sigma^2=\dfrac{E}{(1+\nu)(1-2\nu)}\bigg[(1-\nu)\dfrac{\partial u_2}{\partial y}+\nu\bigg(\dfrac{\partial u_1}{\partial x}+\dfrac{\partial u_3}{\partial z}\bigg) \bigg]\qquad\quad\tau^{13}=\dfrac{E}{2(1+\nu)}\bigg(\dfrac{\partial u_1}{\partial z}+\dfrac{\partial u_3}{\partial x}\bigg)\\
    &\sigma^3=\dfrac{E}{(1+\nu)(1-2\nu)}\bigg[(1-\nu)\dfrac{\partial u_3}{\partial z}+\nu\bigg(\dfrac{\partial u_1}{\partial x}+\dfrac{\partial u_2}{\partial y}\bigg) \bigg]\qquad \quad \tau^{23}=\dfrac{E}{2(1+\nu)}\bigg(\dfrac{\partial u_2}{\partial z}+\dfrac{\partial u_3}{\partial y}\bigg) .
	\end{split}
\end{equation}


The elastic energy related to the internal forces in the configuration corresponding to a generic displacement $\u$ is given by
\begin{equation*}
\mathcal{E}_{el}(\u)=\dfrac{1}{2}\int_\Omega \TT \u: \D \u \, d\x  \, .
\end{equation*}
If we assume that on $\Omega$ act body forces per unit of volume $\f=(f_1,f_2,f_3)$ and boundary forces per unit of surface $\g=(g_1,g_2,g_3)$ we obtain the total energy of the system
\begin{equation} \label{eq:total-energy}
	\mathcal{E}(\u)=\dfrac{1}{2}\int_\Omega \TT \u: \D \u \, d\x-\int_\Omega \f\cdot \u \, d\x-\int_{\partial \Omega} \g\cdot \u \, dS \, .
\end{equation}
Thanks to the symmetry of the stress tensor $\TT \u=(\TT \u)^T$ we infer that for any $\u, \v\in H^1(\Omega;\R^3)$
\begin{equation*}
	\TT \u: \nabla \v =(\TT \u)^T:(\nabla \v)^T=\TT \u:(\nabla \v)^T
\end{equation*}
so that
\begin{equation}\label{eq:Tu-Du}
2(\TT \u: \nabla \v)=\TT \u: \nabla \v+\TT \u:(\nabla \v)^T\quad\Rightarrow\quad   \TT \u: \nabla \v=\TT \u: \D \v \, .
\end{equation}



Recalling the Hooke's law \eqref{Hooke}, we observe that the bilinear form
\begin{equation*}
  (\u,\v) \mapsto \int_\Omega \TT \u:\D \v \, d\x  \, , \qquad  (\u,\v)\in H^1(\Omega;\R^3)\times H^1(\Omega;\R^3)
\end{equation*}
is symmetric, since
\begin{equation} \label{eq:Hooke-2}
  \TT \u:\D \v=\lambda \, (\dive \u) \, (\dive \v)+2\mu \, \D \u:\D \v \, .
\end{equation}
By looking at the total energy $\mathcal E$ in \eqref{eq:total-energy} as a functional $\mathcal E:H^1(\Omega;\R^3)\to \R$ and exploiting the symmetry of the bilinear form above mentioned, we see that a critical point $\u\in H^1(\Omega;\R^3)$ of $\mathcal E$ solves the variational problem
\begin{equation} \label{eq:weak}
  \int_\Omega \TT \u:\D \v \, d\x=\int_\Omega \f\cdot \v \, d\x+\int_{\partial \Omega} \g \cdot \v \, dS \qquad
  \text{for any } \v\in H^1(\Omega;\R^3) \, .
\end{equation}
By \eqref{eq:Tu-Du} and a formal integration by parts, we see that \eqref{eq:weak} is the weak formulation of the boundary value problem
\begin{equation}\label{elastic0}
\begin{cases}
-\dive (\TT \u)=\f\quad&\text{in }\Omega,\\ 
(\TT \u) \n=\g\quad &\text{on }\partial \Omega.
\end{cases}
\end{equation}
Inserting \eqref{Hooke} into \eqref{elastic0} we find the well known equations of linear elasticity
\eqref{pb}.

In the next subsection we prove the existence of solution, stating some classical results about functional spaces of vector valued functions which find a natural application in the theory of linear elasticity.
 These results are related to the well known Korn inequality which has a general validity for vector functions from $\R^N$ to $\R^N$ for any $N \ge 1$. Clearly, in the present paper we will be mainly interested to the case $N=3$, being $\R^3$ the natural space where a solid elastic body can be modelled. For completeness, we will state those results in the general $N$-dimensional case.

\subsection{Existence of a solution}
Let $\Omega\subset \R^n$ a bounded domain, i.e. an open connected bounded set of $\R^N$. Let us introduce the following Sobolev-type space $H^1_\D(\Omega;\R^N)$ defined as the completion of $C^\infty(\overline\Omega;\R^N)$
with respect to the scalar product

\begin{equation}  \label{eq:scal-D}
(\u,\v)_{H^1_\D}=\int_\Omega \D \u:\D \v \, d\x+\int_\Omega \u\cdot \v \, d\x \qquad \text{for any } \u,\v \in  C^\infty(\overline\Omega;\R^N) \, .
\end{equation}

Here $\x=(x_1,\dots,x_N)$ denotes the generic variabile of a function defined in a domain of $\R^N$ and $d\x=dx_1 \dots dx_N$ denotes the $N$-dimensional volume integral in $\R^N$.
From its definition, it is clear that $H^1_\D(\Omega;\R^N)\subset L^2(\Omega;\R^N)$ becomes a Hilbert space with the extension of the scalar product \eqref{eq:scal-D}.

The Korn inequality states the equivalence on the space $C^\infty(\overline\Omega;\R^N)$ between the usual scalar product of $H^1(\Omega;\R^N)$, namely

\begin{equation} \label{eq:scal-H1}
(\u,\v)_{H^1}=\int_\Omega \nabla \u:\nabla \v \, d\x+\int_\Omega \u\cdot \v \, dx \qquad \text{for any } \u,\v \in H^1(\Omega;\R^N) \,
\end{equation}
and the scalar product \eqref{eq:scal-D}. More precisely, if $\Omega \subset \R^N$ is a bounded domain with Lipschitz boundary, then there exists $C>0$ such that
	\begin{equation}\label{t:Korn}
	\int_\Omega |\nabla \u|^2 d\x \le C \left(\int_\Omega |\u|^2 d\x+\int_\Omega |\D \u|^2 d\x\right) \qquad \text{for any } \u\in H^1(\Omega;\R^N) \, .
	\end{equation}

Among the others, for a clear and elegant proof of \eqref{t:Korn}, we address the reader to \cite{korn} by V. A. Kondrat'ev \& O. A. Oleinik.

Thanks to \eqref{t:Korn} we deduce that the Hilbert space $H^1_\D(\Omega;\R^N)$ actually coincides with $H^1(\Omega;\R^N)$ as one can deduce from the definition of $H^1_\D(\Omega;\R^N)$ and the well known result about density of $C^\infty(\overline\Omega;\R^N)$ in $H^1(\Omega;\R^N)$ whenever $\Omega$, is a bounded domain with Lipschitzian boundary.
The Korn inequality is a fundamental tool for proving the existence of weak solutions of \eqref{pb}.

We fix $N=3$ and we observe that by \eqref{eq:Hooke-2} and \eqref{t:Korn}, the bilinear form defined by
\begin{equation}\label{eq:scal-prod-T}
   (\u,\v)_\TT=\int_\Omega \TT \u: \D \v\, d\x+\int_\Omega \u \cdot \v \, d\x \qquad \text{for any } \u,\v \in H^1(\Omega;\R^3)
\end{equation}
is a scalar product in $H^1(\Omega;\R^3)$ which is equivalent to \eqref{eq:scal-H1}.

Let us introduce the space
\begin{equation} \label{eq:V0}
  V_0:=\left\{\v_0 \in H^1(\Omega;\R^3): \int_\Omega \TT \v_0:\D \v \, d\x=0\quad \forall \v\in H^1(\Omega;\R^3)  \right\} \, .
\end{equation}

We observe that $V_0$ coincides with the eigenspace associated to the first eigenvalue of the following eigenvalue problem: $\alpha$ is an eigenvalue if there exists a nontrivial function $\u\in H^1(\Omega;\R^3)$, which will be called eigenfunction associated to $\alpha$, such that
\begin{equation*} 
  \int_\Omega \TT \u:\D \v \, d\x=\alpha \int_\Omega \u\cdot \v \, d\x \qquad \text{for any } \v\in H^1(\Omega;\R^3).
\end{equation*}

In particular if $\alpha=0$ and $\v_0$ is a corresponding eigenfunction we have that

\begin{equation} \label{eq:var-eig}
   \int_\Omega \TT \v_0:\D \v \, d\x=0 \qquad \text{for any } \v\in H^1(\Omega;\R^3) \, .
\end{equation}

After some computation one can verify that $V_0$ is the space of functions $\v=(v_1,v_2,v_3)$ admitting the following representation
\begin{equation} \label{eq:van-cond}
   \begin{cases}
      v_1(x_1,x_2,x_3)=\alpha x_2+\beta x_3+\delta_1  \, ,  \\
      v_2(x_1,x_2,x_3)=-\alpha x_1+\gamma x_3+\delta_2 \, , \\
      v_3(x_1,x_2,x_3)=-\beta x_1-\gamma x_2+\delta_3 \, .
   \end{cases}
\end{equation}
where $\alpha, \beta, \gamma, \delta_1, \delta_2, \delta_3$ are arbitrary constants.

Roughly speaking, configurations associated with such functions $\v\in V_0$ correspond to translations and rotations of the solid body without deforming it in such a way the elastic energy $\mathcal E_{el}(\v)$ equals zero.

Actually, assuming for simplicity $\delta_1=\delta_2=\delta_3=0$, deformations corresponding to displacements $\v\in V_0$ can be considered good approximations of a rotation only for $\alpha, \beta, \gamma$ small; when at least one of the constants $\alpha,\beta, \gamma$ is not small the corresponding deformation of the solid body is no more negligible. In such a case, one may wonder why the elastic energy remains anyway zero; the answer is that in the linear theory only small deformations are allowed so that large deformations are no more meaningful for our model.

Let us recall that we are considering in our model the linearized strain tensor $\D \v$ which is a good approximation of the real strain tensor only for small deformations since the last one also contains quadratic terms in the first order derivatives of $v_1, v_2, v_3$; these quadratic terms can be neglected when first order derivatives are small.

In the next theore we state the existence result for problem \eqref{pb}.
\begin{theorem} \label{p:existence}
	Let $\Omega\subset \R^3$ a bounded domain with Lipschitzian boundary and let $\f\in L^2(\Omega;\R^3)$ and $\g \in L^2(\partial\Omega;\R^3)$. Let us introduce the following compatibility condition
	\begin{equation} \label{eq:compatibility}
	  \int_\Omega \f\cdot \v \, d\x+\int_{\partial\Omega} \g \cdot \v \, dS=0 \qquad \text{for any } \v\in V_0 \, .
	\end{equation}	
Then the following statements hold true:

\begin{itemize}
  \item[$(i)$] problem \eqref{eq:weak} admits at least one solution $\u\in H^1(\Omega;\R^3)$, or equivalently the boundary value problem \eqref{pb} admits at least one weak solution, if and only if \eqref{eq:compatibility} holds true;

  \item[$(ii)$] if $\u\in H^1(\Omega;\R^3)$ is a particular solution of \eqref{eq:weak} then for any $\v_0 \in V_0$ the function $\u+\v_0$ is still a solution of \eqref{eq:weak};

  \item[$(iii)$] if $\u\in H^1(\Omega;\R^3)$ is a particular solution of \eqref{eq:weak} and if $\w\in H^1(\Omega;\R^3)$ is any other solution of \eqref{eq:weak}, then there exists $\v_0 \in V_0$ such that $\w=\u+\v_0$;

  \item[$(iv)$] letting $V_0^\perp$ be the space orthogonal to $V_0$ with respect to the scalar product \eqref{eq:scal-prod-T}, we have that problem \eqref{eq:weak} admits a unique solution in $V_0^\perp$.
\end{itemize}
\end{theorem}

\medskip


\begin{remark}
  The results of Theorem \ref{p:existence} may be extended by replacing the boundary function $\g\in L^2(\partial\Omega;\R^3)$ by an element of the space $H^{-1/2}(\partial \Omega;\R^3)$ denoting the dual space of $H^{1/2}(\partial\Omega;\R^3)$. In such a case, if $\g\in H^{-1/2}(\partial \Omega;\R^3)$ the compatibility condition \eqref{eq:compatibility} has to be replaced by its natural extension
  \begin{equation*} 
	  \int_\Omega \f\cdot \v \, d\x+ \, _{H^{-1/2}} \langle \g, \v\rangle_{H^{1/2}}=0 \qquad \text{for any } \v\in V_0 \, .
	\end{equation*}	
 The validity of this fact comes from the trace theory for vector valued functions ${\bf h}$ admitting a weak divergence in $L^2(\Omega)$; for such functions ${\bf h}$, it is possible to define the trace of ${\bf h}\cdot \n$ as an element of the dual space of $H^{1/2}(\partial\Omega)$, the last one being the space of traces of $H^1(\Omega)$ functions. Such a result has to be applied in our case to each line of $\TT \u$.
\end{remark}

\begin{remark} We observe that, as a consequence of Theorem \ref{p:existence}, if $\u$ and $\w$ are solutions of \eqref{eq:weak}, then the two configurations of the elastic body, corresponding to $\u$ and $\w$, generate the same stress state. More precisely we have $\TT \u=\TT \w$ in $\Omega$ as a consequence of the Hooke's law and of the fact that $\D(\u-\w)$ vanishes in $\Omega$ being $\u-\w\in V_0$. Physically, this is completely reasonable since, given the configuration corresponding to $\u$, the one corresponding to $\w$ can be obtained from the first one by means of rotations and translations of the elastic body, which clearly do not affect the stress state of the solid body itself.
\end{remark}

The proof of Theorem \ref{p:existence} is given in subsection \ref{proofexistence} and is based on standard arguments and the Fredholm alternative.

\section{The hollow cylinder axially loaded at the end faces} 
We consider a circular, finite, homogeneous, isotropic and elastic cylinder with height $h$, radius $b>0$, having a coaxial hole of radius $a>0$.
\begin{figure}[h!]
	\centering
	 \includegraphics[width=5cm]{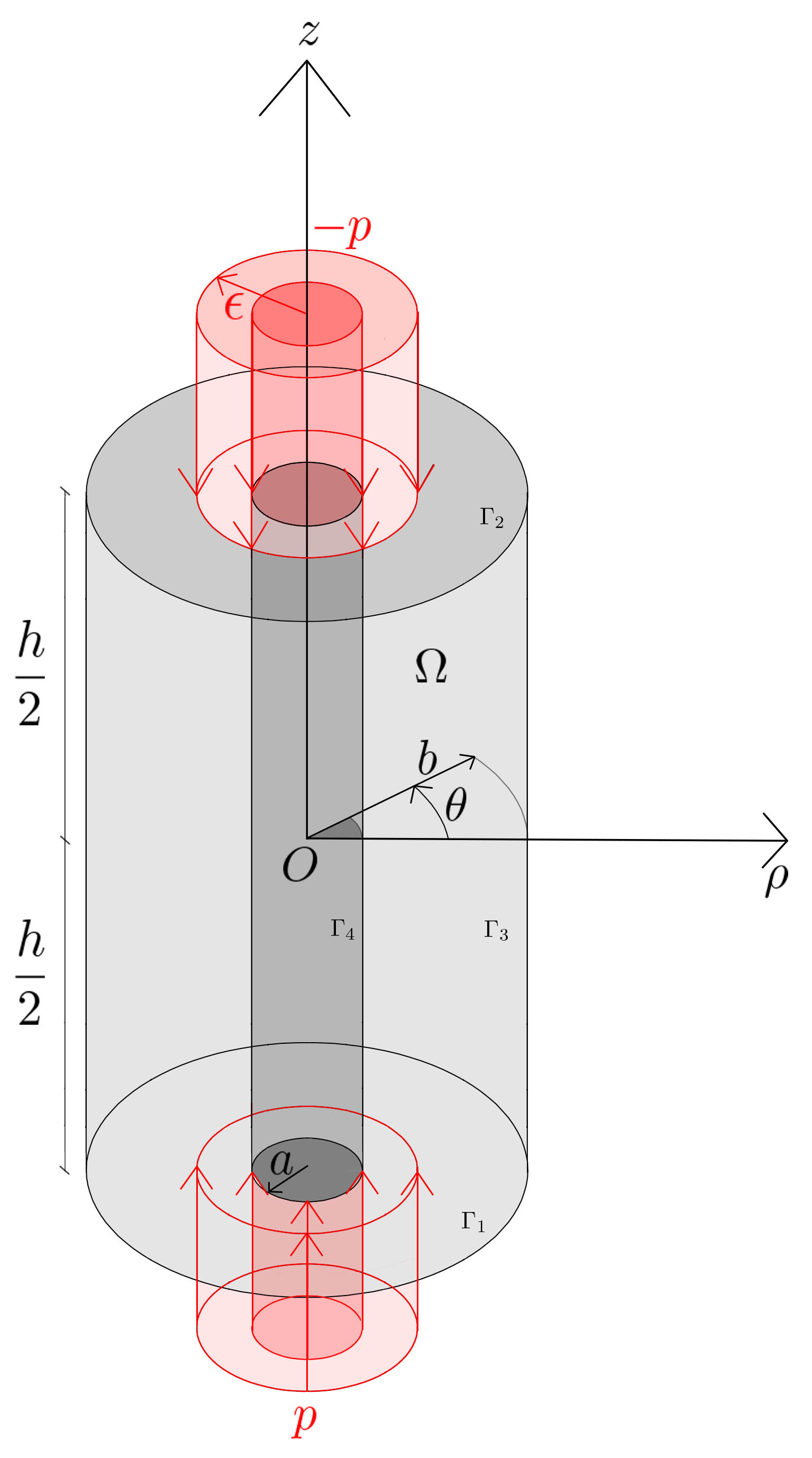}
	\caption{The domain $\Omega$ and in red the load considered. }\label{cylinder0}
\end{figure}
In this section we use the usual notation $x,y,z$ for the three coordinates in $\R^3$. We maintain the notation $d\x$ to denote the differential volume $dxdydz$.

Therefore, we introduce the annular domain $C_{a,b}:=\{(x,y)\in\mathbb{R}^2:a^2<x^2+y^2<b^2\}$ in such a way that
$$
  \Omega=C_{a,b} \times \left(-\tfrac h2,\tfrac h2\right) \, .
$$

In the sequel we want to model a hollow cylinder subject to an external load acting on the upper and lower faces of the cylinder compressing the cylinder itself. Recalling the notations introduced in \eqref{pb}, we will then assume that the volume forces represented by the vector function $\f$ vanish everywhere in $\Omega$.

In order to better describe the surface forces represented by the vector function $\g$, we split $\partial\Omega$ in four regular parts
$$
   \Gamma_1:=C_{a,b} \times \left\{-\tfrac h 2 \right\} \, , \qquad \Gamma_2:=C_{a,b} \times \left\{\tfrac h 2 \right\} \, ,
$$
$$
 \Gamma_3:=\{(x,y)\in\mathbb{R}^2:x^2+y^2=b^2\}\times \left(-\tfrac h 2,\tfrac h 2 \right) \, , \qquad \Gamma_4:=\{(x,y)\in\mathbb{R}^2:x^2+y^2=a^2\}\times \left(-\tfrac h 2,\tfrac h 2\right),
$$
having respectively outward unit normal vectors $(0,0,-1)$, $(0,0,1)$, $(x/b,y/b,0)$ when $(x,y,z)\in \Gamma_3$ and $(-x/a,-y/a,0)$ when $(x,y,z)\in \Gamma_4$.
In this way, the outward unit normal vector $\n$ is well defined on the whole $\partial\Omega$.

Exploiting the above notations, the vector function $\g$ can be represented in the following way
\begin{equation} \label{eq:def-g}
  \g(x,y,z)=
    \begin{cases}
       (0,0,\chi_p(x,y))  & \text{for any }  (x,y,z)\in \Gamma_1,  \\
       (0,0,-\chi_p(x,y))   & \text{for any }  (x,y,z)\in \Gamma_2,  \\
       (0,0,0)           & \text{for any }  (x,y,z)\in \Gamma_3\cup \Gamma_4,
    \end{cases}
\end{equation}
where the function $\chi_p:C_{a,b} \to \R$, $p\in \R_+$, is defined by
\begin{equation}\label{load0}
\chi_p(x,y):=
  \begin{cases}
     p  & \qquad \text{if } a^2\leq x^2+y^2<\epsilon^2   , \\
     0  & \qquad \text{if } \epsilon^2<x^2+y^2 \leq b^2  ,
\end{cases}
\end{equation}
for some $\epsilon \in (a,b)$.

Resuming all the assumptions on $\f$ and $\g$ we are led to consider the problem
\begin{equation} \label{eq:weak-2}
  \begin{cases}
     -\mu \Delta \u-(\lambda+\mu)\nabla (\dive \u)=\0  & \qquad \text{in } \Omega \, , \\[3pt]
     (\TT \u)\n =(0,0,\chi_p)                           & \qquad \text{on } \Gamma_1, \\[3pt]
     (\TT \u)\n =(0,0,-\chi_p)                           & \qquad \text{on } \Gamma_2, \\[3pt]
         (\TT \u)\n = \0                                   & \qquad \text{on }  \Gamma_3\cup \Gamma_4 \, . 
  \end{cases}
\end{equation}

Among all solutions of \eqref{eq:weak-2} which can be obtained by a single solution by adding to it a function in the space $V_0$, we focus our attention on the unique solution $\u=(u_1,u_2,u_3)$ of \eqref{eq:weak-2} in the space $V_0^\perp$ where orthogonality is meant in the sense of the scalar product defined in \eqref{eq:scal-prod-T}, see Theorem \ref{p:existence}.
From a geometric point of view, condition $\u \in V_0^\perp$ avoids translations and rotations of the hollow cylinder, being $V_0$ the space of displacement functions which generate translations and rotations.

In subsection \ref{proofprop0} we prove a symmetry result for the unique solution of \eqref{eq:weak-2} in the space $V_0^\perp$, whose validity is physically evident, but which however needs a rigorous proof:

\begin{proposition}\label{prop0}
   Let $\u$ be the unique solution of \eqref{eq:weak-2} in the space $V_0^\perp$. Then $\u$ satisfies the following symmetry properties:

   \begin{itemize}
     \item[$(i)$] for any $(x,y,z)\in \Omega$ we have
	\begin{equation}\label{eq:symmetry-1}
	  u_1(x,y,-z)=u_1(x,y,z), \quad u_2(x,y,-z)=u_2(x,y,z) \quad u_3(x,y,-z)=-u_3(x,y,z) \, ,
	\end{equation}
\begin{equation}\label{eq:symmetry-2}
	  u_1(-x,y,z)=-u_1(x,y,z), \quad u_2(-x,y,z)=u_2(x,y,z) \quad u_3(-x,y,z)=u_3(x,y,z) \, ,
\end{equation}
\begin{equation}\label{eq:symmetry-3}
	  u_1(x,-y,z)=u_1(x,y,z), \quad u_2(x,-y,z)=-u_2(x,y,z) \quad u_3(x,-y,z)=u_3(x,y,z) \, ;
	\end{equation}

\item[$(ii)$] the third component $u_3$ of the solution $\u$ is axially symmetric in the sense that:
\begin{equation*} 
    u_3(x_1,y_1,z)=u_3(x_2,y_2,z) \quad \text{$\forall$ $(x_1,y_1,z), (x_2,y_2,z)\in \Omega$ with $x_1^2+y_1^2=x_2^2+y_2^2$} \, ;
\end{equation*}

\item[$(iii)$] the first two components $u_1$, $u_2$ of the solution $\u$ form a central vector field in two dimensions
in the sense that
\begin{equation*}
   \left|(u_1,u_2)\right|_{|(x_1,y_1,z)}=\left|(u_1,u_2)\right|_{|(x_2,y_2,z)} \quad \text{$\forall$ $(x_1,y_1,z), (x_2,y_2,z)\in \Omega$ with $x_1^2+y_1^2=x_2^2+y_2^2$}
\end{equation*}
and
\begin{equation*}
   (u_1,u_2)_{|(x,y,z)}=\left|(u_1,u_2)\right|_{|(x,y,z)} \, \left(\frac{x}{\sqrt{x^2+y^2}},\frac{y}{\sqrt{x^2+y^2}}\right) \quad \text{for any } (x,y,z)\in \Omega \, .
\end{equation*}

 \end{itemize}
\end{proposition}

\subsection{Periodic extension of the problem}\label{periodic}


Our next purpose is to look for and construct a solution $\u=(u_1,u_2,u_3)$ of \eqref{eq:weak-2} admitting a Fourier series expansion and, hence, admitting a periodic extension defined on the whole $C_{a,b} \times \R$.
In order to obtain this construction, we need to assume that the horizontal displacements $u_1$ and $u_2$ vanish on the upper and lower faces of the hollow cylinder $\Omega$:
\begin{equation} \label{eq:vanishing-1-2}
    u_1\left(x,y,\tfrac h2\right)=u_1\left(x,y,-\tfrac h2\right)=0 \quad \text{and}
    \quad u_2\left(x,y,\tfrac h2\right)=u_2\left(x,y,-\tfrac h2\right)=0 \, .
\end{equation}
We find a solution satisfying \eqref{eq:vanishing-1-2} and, a posteriori, we show that it necessarily coincides with the unique solution of \eqref{eq:weak-2} belonging to $V_0^\perp$, see the end of the proof of Theorem \ref{teo_soluzioni}.

As a first step, since $\u\in H^1(\Omega;\R^3)$, we define a function, still denoted for simplicity by $\u$, on the domain $C_{a,b} \times \left(-\frac h 2,\frac{3h}{2}\right)$ by extending it in suitable way: the new function $\u$ coincides with the original function $\u$ on $C_{a,b} \times \left(-\frac h 2,\frac{h}{2}\right)$ and
\begin{equation} \label{eq:sym-ext-123}
  u_1(x,y,z)=-u_1(x,y,h-z) \, , \quad u_2(x,y,z)=-u_2(x,y,h-z) \, , \quad u_2(x,y,z)=u_3(x,y,h-z)  \, ,
\end{equation}
for any $(x,y,z)\in C_{a,b} \times \left(\tfrac h 2, \tfrac{3h}{2}\right)$.
This means that $u_1$ and $u_2$ are antisymmetric with respect to $z=\frac h2$ and $u_3$ is symmetric with respect to $z=\frac h2$.
This symmetric extension with respect to $z=\frac h 2$ produces a function $\u \in H^1\left(C_{a,b} \times \left(-\frac h 2, \frac{3h}{2} \right);\R^3\right)$ thanks to condition \eqref{eq:vanishing-1-2}.

The second step is to extend the new function $\u: C_{a,b} \times \left(-\frac h 2,\frac{3h}{2}\right)\to \R$ to the whole $C_{a,b} \times \R$ as a $2h$-periodic function in the variable $z$.
It is easy to understand that the periodic extension, still denoted for simplicity by $\u$, is a function satisfying $\u \in H^1(C_{a,b} \times I;\R^3)$ for any open bounded interval $I$.

The periodic extension of the boundary data can be achieved according to the next lemma, proved in subsection \ref{prooflemma1}. We state here some lemmas in order to understand the main steps in the construction of the solution of \eqref{eq:weak-2}, given in the final theorem.
\begin{lemma}\label{lemma1}
Let $\u$ be the periodic extension of the solution of \eqref{eq:weak-2} defined as above and let ${\bm \Lambda}$ be the distribution defined by 
	\begin{equation} \label{eq:distr-equation}
		-{\rm div}(\TT  \u)={\bm \Lambda}     \qquad \text{in } \mathcal D'(C_{a,b}\times \R;\R^3) \, .
	\end{equation}
Then ${\bm \Lambda}$ admits the following Fourier series expansion
	\begin{equation} \label{eq:distr-expansion}
		{\bm \Lambda}=(\Lambda_1,\Lambda_2,\Lambda_3)=\left(0,0,\chi_p(x,y) \,  \sum_{m=0}^{+\infty} (-1)^{m+1} \, \frac{4}{h} \sin\left[\frac \pi h (2m+1)z\right] \right).
	\end{equation}
\end{lemma}

The symmetry properties \eqref{eq:sym-ext-123} and the construction of the periodic extension, allow expanding $\u=(u_1, u_2,  u_3)$ in Fourier series with respect to the variable $z$: 

\begin{align} \label{eq:comp-expansion}
	&	 u_1(x,y,z)=\sum_{k=0}^{+\infty} \varphi_k^1(x,y)\cos\left(\tfrac \pi h \, kz\right) \, , \qquad \quad u_2(x,y,z)=\sum_{k=0}^{+\infty} \varphi_k^2(x,y)\cos\left(\tfrac \pi h \, kz\right) \, , \\[7pt]
	\notag &   u_3(x,y,z)=\sum_{k=0}^{+\infty} \varphi_k^3(x,y)\sin\left(\tfrac \pi h \, kz\right) \, .
\end{align}

 In subsection \ref{prooflemma2} we prove the following lemma.
 \begin{lemma}\label{lemma2}
For any $k\ge 1$ odd, there exists a unique
$(\varphi_k^1,\varphi_k^2,\varphi_k^3)\in H^1( C_{a,b};\R^3)$, satisfying \eqref{eq:weak-2} and \eqref{eq:comp-expansion}. For any $k\ge 2$ even, there exists a unique trivial $(\varphi_k^1,\varphi_k^2,\varphi_k^3) \equiv (0,0,0)$ in $C_{a,b}$, satisfying \eqref{eq:weak-2} and \eqref{eq:comp-expansion}.
 \end{lemma}

\begin{remark}
We observe that for $k=0$, the boundary value problem \eqref{eq:weak-2}, or equivalently \eqref{eq:2D}-\eqref{eq:BC-2D} and \eqref{eq:big-var-2}, see the proof in subsection \ref{prooflemma2}, admits an infinite number of solutions. 
More precisely, these solutions are in form $(\varphi_0^1,\varphi_0^2,\varphi_0^3)=(c_1,c_2,c_3)$ where $c_1,c_2,c_3$ are three arbitrary constants. We may choose $\varphi_0^3\equiv 0$ being irrelevant in the Fourier expansion of $u_3$.
Concerning the other two components, we have necessarily $\varphi_0^1\equiv \varphi_0^2\equiv 0$ in $C_{a,b}$ due to the odd symmetry of $u_1$ and $u_2$ with respect to the variables $x$ and $y$, as stated in \eqref{eq:symmetry-2} and \eqref{eq:symmetry-3}.
\end{remark}

\subsection{Cylindrical coordinates exchange}\label{cil}
The symmetry properties of $\u$ stated in Proposition \ref{prop0} imply that $\varphi_k^3$ is a radial function and the vector field $(\varphi_k^1,\varphi_k^2)$ is a central vector field in the plane, in the sense that it is oriented toward the origin and its modulus is a function only of the distance from the origin. This implies that for any $k\ge 1$ odd, there exist two radial functions $Y_k=Y_k(\rho)$ and $Z_k=Z_k(\rho)$ such that in polar coordinates we may write
\begin{equation} \label{eq:Yk-Zk}
	\varphi_k^1(\rho,\theta)=Y_k(\rho)\cos\theta \, , \qquad
	\varphi_k^2(\rho,\theta)=Y_k(\rho)\sin\theta \, , \qquad
	\varphi_k^3(\rho,\theta)=Z_k(\rho)  \, ,
\end{equation}
with $\rho \in[a,b]$ and $\theta \in [0,2\pi)$.

In Section \ref{proofs} we show that $Y_k$ and $Z_k$ solve a proper boundary value problem. More precisely, this fact will be shown in subsection \ref{ss:p-33} which is devoted to the proof of the next lemma, where we state existence and uniqueness for solutions of the boundary value problem mentioned above.

\begin{lemma} \label{p:prop-ex-un-rad}
  Let $\Psi_k: C_{a,b}\to \R$ be defined as
  \begin{equation} \label{eq:Psi-k}
  	\Psi_k(x,y):=
  	\begin{cases}
  		\ds{(-1)^{\frac{k+1}{2}} \frac{4}{h} \,  \chi_p(x,y) } & \qquad \text{if $k$ is odd,}  \\[5pt]
  		0                                                       & \qquad \text{if $k$ is even,}
  	\end{cases}\qquad\forall(x,y)\in  C_{a,b}.
  \end{equation}
  For any $k\ge 1$ odd, the boundary value problem
  \begin{equation}\label{pbcil3}
  \begin{cases}
  Y_k''(\rho)+\dfrac{Y_k'(\rho)}{\rho}-\dfrac{Y_k(\rho)}{\rho^2}-\dfrac{\mu}{\lambda+2\mu} \, \dfrac{\pi^2 k^2}{h^2} \, Y_k(\rho)+\dfrac{\lambda+\mu}{\lambda+2\mu} \, \dfrac{\pi k}h Z_k'(\rho)=0 \quad & \text{in } ( a, b) \, , \\[12pt]
  Z_k''(\rho)+\dfrac{Z_k'(\rho)}{\rho}-\dfrac{\lambda+2\mu}{\mu} \, \dfrac{\pi^2 k^2}{h^2} Z_k(\rho)-\dfrac{\lambda+\mu}{\mu} \, \dfrac{\pi k}{h} \bigg[Y_k'(\rho)+\dfrac{Y_k(\rho)}{\rho}\bigg]=-\dfrac 1 \mu \, \Psi_k(\rho)  \quad & \text{in }( a, b) \, , \\[12pt]
  (\lambda+2\mu)Y_k'(\rho)+\frac{\lambda}{\rho} \, Y_k(\rho)+\lambda \frac{\pi k}{h} Z_k(\rho)=0  \, , \quad& \rho\in \{ a, b\} \\[12pt]
  Z_k'(\rho)-\frac{\pi k}{h} Y_k(\rho)=0 \, , \quad& \rho\in \{ a, b\}
  \end{cases}
  \end{equation}
  admits a unique solution $(Y_k,Z_k)\in H^1( a, b;\R^2)$.
\end{lemma}

About existence and uniqueness of solutions of \eqref{pbcil3}, in subsection \ref{ss:p-33} we only give an idea of the proof since it can be proved exactly as Lemma \ref{lemma2} of which
Lemma \ref{p:prop-ex-un-rad} is the radial version.

Now we need a more explicit representation for the unique solution $(Y_k,Z_k)$ of \eqref{pbcil3}. This will be done by performing a power series expansion in which the coefficients will be characterized explicitly in terms of a suitable iterative scheme. As a byproduct of this result in Section \ref{num} we also obtain a numerical approximation of the exact solution and we estimate the corresponding error. Being a linear problem, we proceed by applying the superposition principle and we provide the explicit formula in the next lemma.

\begin{lemma}\label{lemma3}
	For any $k\geq 1$, odd, let ${\bm \Upsilon}_k=(Y_k,Z_k)$ the unique solution of \eqref{pbcil3}. Omitting for brevity the $k$-index, we have a unique $(C_1,C_2,C_3,C_4)\in \R^4$ such that  
	\begin{equation}\label{sol1}
	{\bm \Upsilon}(\rho)=C_1 {\bm \Upsilon}^1(\rho)+C_2 {\bm \Upsilon}^2(\rho)+C_3 {\bm \Upsilon}^3(\rho)+C_4 {\bm \Upsilon}^4(\rho)+\overline{{\bm \Upsilon}}(\rho),
	\end{equation}
	where ${\bm \Upsilon}^j=(Y^j,Z^j)$ with $j=1,\dots,4$ are four linear independent solutions of the corresponding homogenous system and $\overline{{\bm \Upsilon}}=(\overline Y,\overline Z)$ solves
	\begin{equation} \label{eq:Upsilon-bar}
	\begin{pmatrix}
	\overline Y(\rho) & 
	\overline Y'(\rho)& 
	\overline Z(\rho) & 
	\overline Z'(\rho)
	\end{pmatrix}^T
	={\bf W}(\rho) \int_{ a}^\rho ({\bf W}(r))^{-1}
	\begin{pmatrix}
	0 & 
	0 & 
	0 & 
	-\frac 1\mu \Psi_k(r)
	\end{pmatrix}^T \, dr \, , \qquad \rho>0 \, ,
	\end{equation}
	being ${\bf W}(\rho)$ the wronskian obtained through ${\bm \Upsilon}^j(\rho)$ ($j=1,\dots,4$).
	Each of the linear independent solutions of the homogeneous system can be written as
	\begin{equation}\label{expansion2}
	\begin{cases}
	{\ds Y^j(\rho)=\sum_{n=-1}^{+\infty} a^j_n \, \rho^n+(\ln \rho) \sum_{n=0}^{+\infty} b^j_n \, \rho^n} \, , \\[15pt]
	{\ds Z^j(\rho)=\sum_{n=0}^{+\infty} c^j_n \, \rho^n+(\ln \rho) \sum_{n=0}^{+\infty} d^j_n \, \rho^n} \,
	\end{cases}\qquad (j=1,\dots,4),
	\end{equation}
	where the coefficients are uniquely determined.
\end{lemma}

In the proof of the Lemma we give all the details related to the computation of the constants $C_j$ in \eqref{sol1} and of the coefficients in the series \eqref{expansion2}, see subsection \ref{prooflemma3}. As a consequence of Lemmas \ref{lemma1}-\ref{lemma2}-\ref{p:prop-ex-un-rad}-\ref{lemma3} we state the main theorem, whose proof can be found in subsection \ref{ss:main-teo}.

\begin{theorem}\label{teo_soluzioni}
  Let $\u$ be the unique solution of \eqref{eq:weak-2} satisfying $\u \in V_0^\perp$ and let $(Y_k,Z_k)$, $k \ge 1$ odd, be the unique solution of \eqref{pbcil3}. Then, in cylindrical coordinates, $\u=(u_1,u_2,u_3)$ admits the following representation:
  \begin{equation} \label{eq:series-expansions}
    \begin{cases}
      \ds{u_1(\rho,\theta,z)= \sum_{m=0}^{+\infty} Y_{2m+1} \left(\rho\right) \cos \theta \cos\left[\tfrac{(2m+1)\pi}h \, z\right]  \, , }  \\[15pt]
      \ds{u_2(\rho,\theta,z)= \sum_{m=0}^{+\infty} Y_{2m+1} \left(\rho\right) \sin \theta \cos\left[\tfrac{(2m+1)\pi}h \, z\right]  \, , }  \\[15pt]
      \ds{u_3(\rho,\theta,z)= \sum_{m=0}^{+\infty} Z_{2m+1} \left(\rho\right) \sin\left[\tfrac{(2m+1)\pi}h \, z\right]  \, ,             }
    \end{cases}
  \end{equation}
  with $\rho\in (a,b)$, $\theta\in [0,2\pi)$, $z\in \left(-\frac h 2, \frac h 2\right)$ where the three series in \eqref{eq:series-expansions} converge weakly in $H^1(\Omega)$ and strongly in $L^2(\Omega)$.

  Moreover, letting ${\bm U}_M=(U_M^1,U_M^2,U_M^3)$ be the sequence of vector partial sums corresponding to the series expansions in \eqref{eq:series-expansions}, we have for any $M\ge 1$
  \begin{align} \label{eq:stima-L2}
    & \|U_M^1-u_1\|_{L^2(\Omega)}\le \frac{pb^2}{\mu} \sqrt{\frac{h (b-a)}{2a\pi}} \, \frac 1{\sqrt M}  \, , \qquad
    \|U_M^2-u_2\|_{L^2(\Omega)}\le \frac{pb^2}{\mu} \sqrt{\frac{h (b-a)}{2a\pi}} \, \frac 1{\sqrt M} \, , \\[5pt]
    \notag & \|U_M^3-u_3\|_{L^2(\Omega)}\le \frac{p}{\mu a\pi^2 } \sqrt{ \frac{h^3 b^3 (b-a)}{24}} \, \frac{1}{\sqrt{M^3}} \, .
  \end{align}

\end{theorem}

\section{An engineering application}\label{num}

\begin{figure}[h]
	\centering
	\includegraphics[width=10cm]{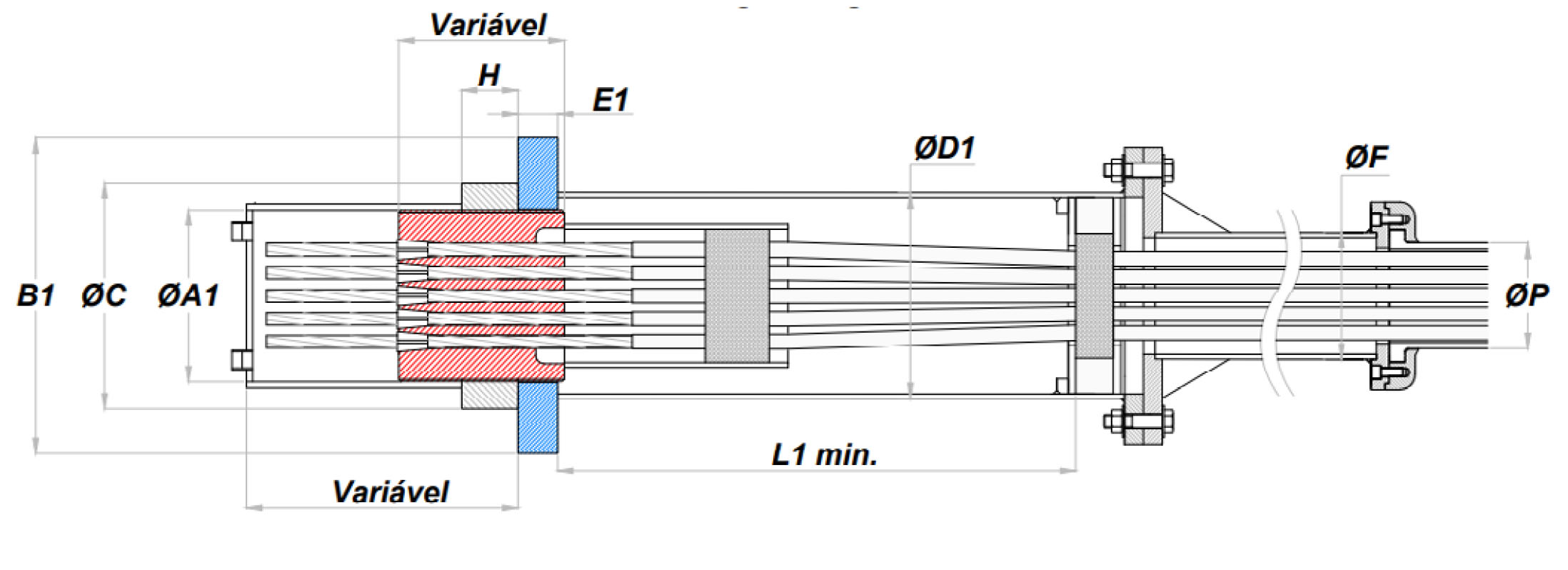}\quad \includegraphics[width=10cm]{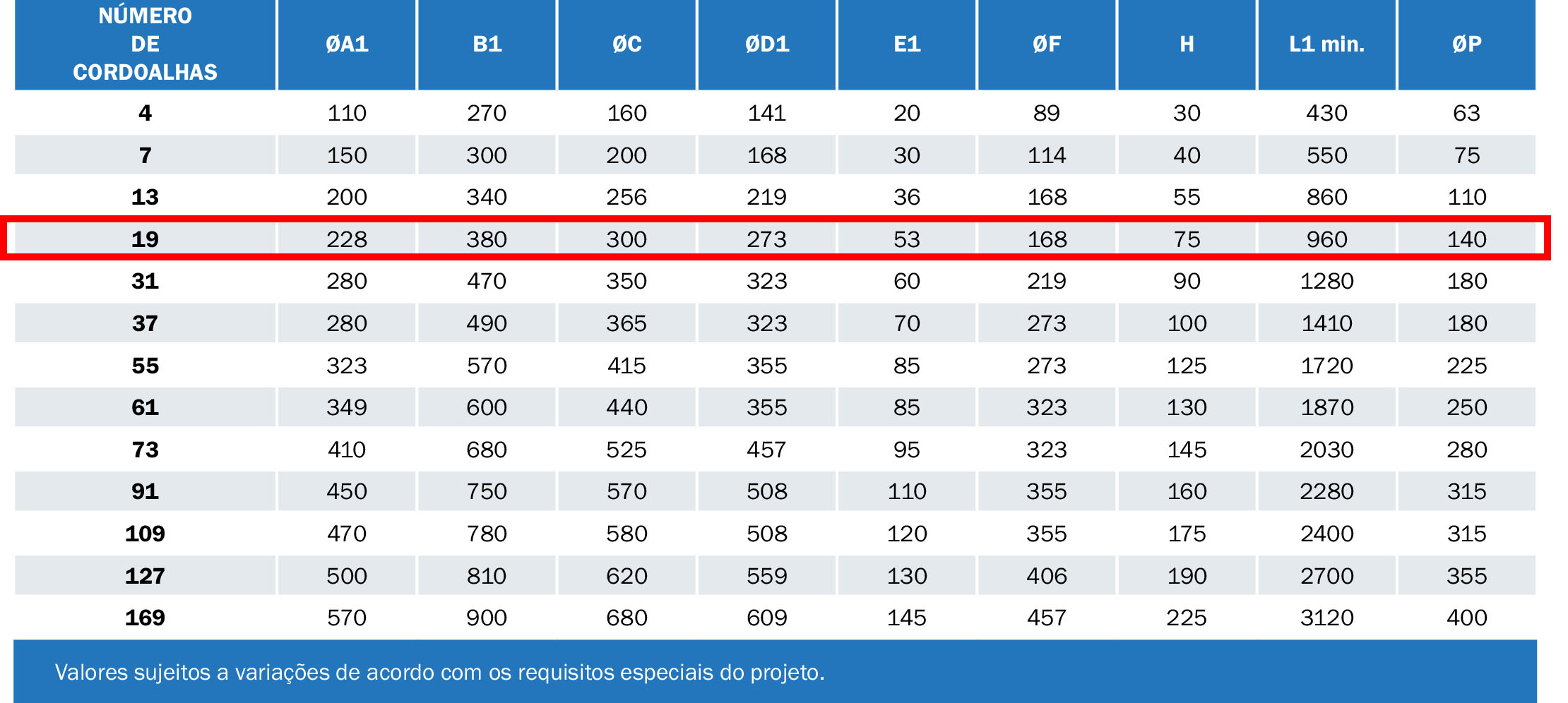}
	\caption{Detail of the strands anchorage and in table the geometric features for a 19 strands element, from the commercial catalogue \cite{protende}.}\label{fig:detail}
\end{figure}

In this section we consider a case of study: a hollow cylinder having the features of a blister for the bridge in Figure \ref{fig:blister}. In Table \ref{tab1} we give the mechanical parameters, see also Figure \ref{cylinder0}.
We consider stays composed of 19 strands, see Figure \ref{fig:detail}, suitable to bear the concentrated load $P$ in Table \ref{tab1}.
$P$ is computed from the executive project, while the diameter $2a$ is taken from the catalogue of Protende ABS-2021 \cite{protende}, a company producing such elements, see in Figure \ref{fig:detail} the diameter $\phi D_1$ for 19 strands anchorage; 
\begin{table}[h!]\centering
	\begin{tabular}{ccl}
		\hline
		$h$ & 3.00 $\,\,$ m & Height of the cylinder\\
		$2a$ & 273 mm& Diameter of the cylindrical hollow\\
		$2b$ & 800 mm & External diameter of the cylinder\\
		$2\epsilon$ & 425 mm & External diameter of the load\\
		$P$ & 1900 kN & Concentrated load\\
		$E$ & 35000 MPa & Young modulus of the concrete\\
		$\nu$ & 0.2 & Poisson ratio of the concrete\\
		\hline
	\end{tabular}
	\caption{Mechanical parameters assumed.}\label{tab1}
\end{table}
hence, the distributed load in \eqref{load0} is given by $p=\frac{P}{\pi(\epsilon^2-a^2)}=22.80$ MPa.

Our purpose is to obtain a good approximation of the functions ${\bm \Upsilon}^j=(Y^j,Z^j)$, $j\in \{1,2,3,4\}$ introduced in Lemma \ref{lemma3}. For $j=1,\dots,4$ we consider the approximate solution ($N\geq1$)
\begin{equation}\label{approx}
Y_{N}^j(\rho)=\sum_{n=-1}^{N}  a_n \, \rho^n+(\ln \rho) \sum_{n=0}^{N}  b_n \, \rho^n  \quad \text{and} \quad
Z_{N}^j(\rho)=\sum_{n=0}^{N-1}  c_n \, \rho^n+(\ln \rho) \sum_{n=0}^{N-1}  d_n \, \rho^n \, .
\end{equation}

The reason for in \eqref{approx} we have $n=0,\dots,N-1$ in the expansion of $Z_N^j$ will be clarified in the proof in subsection \ref{prooferror}
of the next proposition about an estimate of the truncating error.
\begin{proposition}\label{p:error}
	Let $k>1$ , $k \in \mathbb{N}$ odd, and let $N\geq 3$, odd integer, be the truncating index of the series as in \eqref{approx}. Then, letting
\begin{equation*}
   E_{k,N}:=\max_{j\in \{1,2,3,4\}} \left\{ \max \left\{\max_{\rho \in [a,b]} |Y^j_{N}(\rho)-Y^j(\rho)|, \max_{\rho \in [a,b]} |Z^j_{N}(\rho)-Z^j(\rho)| \right\} \right\} \, ,
\end{equation*}
we have that
	\begin{equation}\label{err}
	E_{k,N}\leq  C(a,b,k) \dfrac{3(2\lambda+5\mu)(\lambda+\mu)^2}{16\mu^3} \, \left(\frac{\pi kb}h \right)^{N+2}
e^{\left(\frac{\pi kb}h \right)^2}\, \frac{ (N+3)(3N^3+21N^2+42N +32)}{2^{N} \left[\left(\frac{N+1}{2}\right)!\right]^2},
	\end{equation}
	where
$$
  C(a,b,k)=\max\left\{1,\tfrac{h}{\pi kb}\right\}
\max\Big\{1,\left|\ln\left(\tfrac{\pi k a}{h} \right) \right|,
|\ln\left(\tfrac{\pi k b}{h} \right)|\Big\}
\max\Big\{\tfrac{\pi k}h ,\tfrac{\mu}{\lambda+\mu}\, \tfrac{\pi k}h \ln\left(\tfrac{\pi k}h \right), \tfrac{2(\lambda+2\mu)}{\lambda+\mu} \, \tfrac h {\pi k} \ln\left(\tfrac{\pi k}h \right) \Big\} \, .
$$
\end{proposition}

Once we have \eqref{err}, one may choose $N$ in such a way that
\begin{equation} \label{eq:condition-trunc}
\frac{\ds{ E_{k,N}  }}{\ds{\min_{j\in\{1,2,3,4\}}\left\{\min \left\{\max_{\rho\in [a,b]} |Y_{N}^j(\rho)|, \ \max_{t\in [a,b]} |Z_{N}^j(\rho)| \right\}\right\}}}<\varepsilon
\end{equation}
with $\varepsilon$ small enough. Condition \eqref{eq:condition-trunc} means that the truncation error is relatively small compared to the order of magnitude of both functions $ Y^j_{N}$ and $Z^j_{N}$ for all $j\in \{1,2,3,4\}$.

\begin{figure}[h!]
	\centering
	  \includegraphics[width=17.5cm]{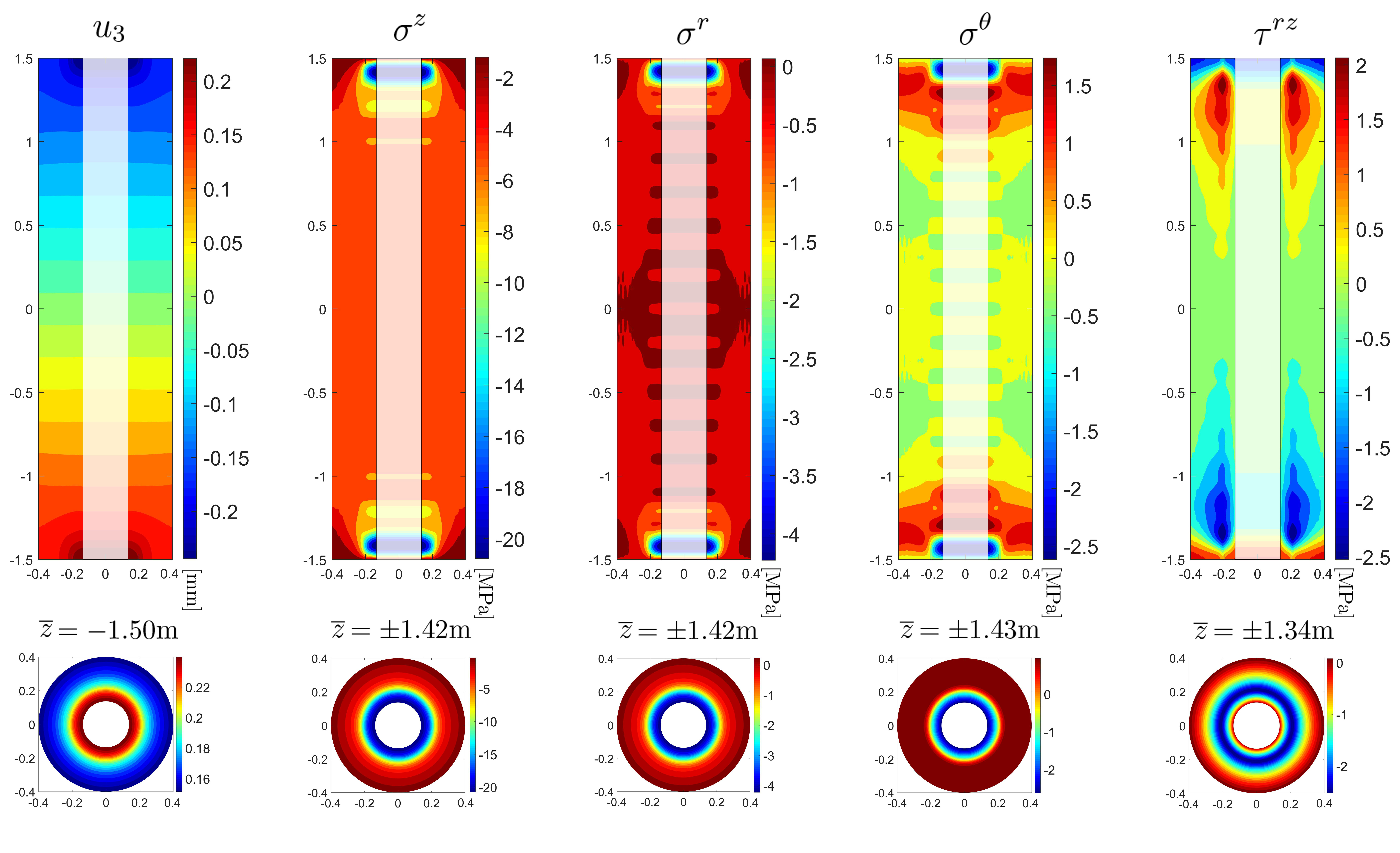}
	\caption{From the the left the vertical displacement $u_3$ in mm, the vertical stress $\sigma^z$ in MPa, the radial stress $\sigma^r$ in MPa, the angular stress $\sigma^\theta$ in MPa and  the tangential stress $\tau^{rz}$ in MPa.}\label{plots}
\end{figure}

In our numerical simulation the condition \eqref{eq:condition-trunc} is verified by making use of estimate \eqref{err} on the truncation error $E_{k,N}$, i.e. the program verifies at each step the validity of  \eqref{eq:condition-trunc} in which the numerator of the fraction  is replaced by the majorant in \eqref{err}. The program runs until the value of $N$ is sufficiently large to guarantee \eqref{eq:condition-trunc}.

In Figure \ref{plots} we plot a vertical section of the cylinder and the corresponding more stressed horizontal section. We show the vertical displacement $u_3$ and the following components of the stress tensor in cylindrical coordinates
\begin{equation}\label{stress2}
\begin{split}
&\sigma^z=\dfrac{2\mu}{1-2\nu}\bigg[(1-\nu)\dfrac{\partial u_{3}}{\partial z}+\nu \bigg(\frac{u^r}{\rho}+\dfrac{\partial u^r}{\partial \rho}\bigg)\bigg]\qquad \sigma^r=\dfrac{2\mu}{1-2\nu}\bigg[(1-\nu)\dfrac{\partial u^r}{\partial \rho}+\nu \bigg(\frac{u^r}{\rho}+\dfrac{\partial u_3}{\partial z}\bigg)\bigg]\\
&\sigma^\theta=\dfrac{2\mu}{1-2\nu}\bigg[(1-\nu)\frac{u^r}{\rho}+\nu \bigg(\dfrac{\partial u^r}{\partial \rho}+\dfrac{\partial u_3}{\partial z}\bigg)\bigg]\qquad \tau^{r z}=\mu \bigg[\dfrac{\partial u^r}{\partial z}+\dfrac{\partial u_3}{\partial \rho}\bigg],
\end{split}
\end{equation}
where $u^r=\sqrt{u_1^2+u_2^2}$ is the radial displacement.
We point out that putting ${\bf n}=(\cos \theta,\sin \theta, 0)$, ${\bf t}=(-\sin \theta, \cos \theta,0)$ and ${\bf k}=(0,0,1)$, the four components introduced in \eqref{stress2} are defined by $\sigma^z:=(\TT \u){\bf k}\cdot {\bf k}$, $\sigma^r:=(\TT \u){\bf n}\cdot {\bf n}$, $\sigma^\theta:=(\TT \u){\bf t}\cdot {\bf t}$ and
$\tau^{rz}:=(\TT \u){\bf n}\cdot {\bf k}$ and the representation \eqref{stress2} can be deduced by \eqref{eq:stress} and \eqref{stresses}.

We consider an approximate solution ${\bm U}_M$ as stated in Theorem \ref{teo_soluzioni} truncating the Fourier series at $M=29$ with $\eps<10^{-3}$ in \eqref{eq:condition-trunc}, implying $N=123$ in \eqref{approx} and $\|U_{29}^1-u_1\|_{L^2(\Omega)}\leq 4.46\cdot 10^{-5}$ m$^{5/2}$, $\|U_{29}^3-u_3\|_{L^2(\Omega)}\leq 1.02\cdot 10^{-6}$ m$^{5/2}$ in \eqref{eq:stima-L2}. In Table \ref{tab} we give the maximum absolute values of the variables involved, including the coordinate of the point $(\overline \rho, \overline z)$ where they are assumed (for all $\theta\in [0,2\pi)$ thanks to the radial symmetry of the problem).
\begin{table}[h]\centering
	\begin{tabular}{cccc}
		\hline
		&$\max|\cdot|$&$\overline \rho$&$\overline z$ [m]\\
		\hline
		$u_3$ & $0.24$ mm & $a$ & $\pm$ 1.50\\
		$\sigma^z$ & $20.77$ MPa & $a$ & $\pm$ 1.42\\
		$\sigma^r$ & $4.23$ MPa & $a$ & $\pm$1.42\\
		$\sigma^{\theta}$ & $2.62$ MPa & $a$ & $\pm$ 1.43\\
		$\tau^{rz}$ & $2.52$ MPa & $\epsilon$ & $\pm$ 1.34\\
		\hline
	\end{tabular}
	\caption{Maximum absolute values and points of $\Omega$ in which they are assumed.}\label{tab}
\end{table}

As expected the vertical displacement $u_3$ achieves its maximum absolute value at $z=\pm \frac{h}{2}$. From the plots we see that there are two (symmetric) critical zones where we observe the loading diffusion; they are close to the upper and bottom faces of the cylinder and involve approximately the 20\% of the closest volume, i.e. the volume of $\Omega$ such that $z\in(-\frac{h}{2},-\frac{2h}{5})\cup(\frac{2h}{5},\frac{h}{2}) $.

\section{Proofs of the results}\label{proofs}

\subsection{Proof of Theorem \ref{p:existence}}\label{proofexistence}

By identity \eqref{eq:Hooke-2}, the estimates $(\dive \u)^2 \le |\nabla \u|^2$ and $|\D \u|^2\le |\nabla \u|^2$ and the H\"older inequality we infer that for any $\u, \v \in H^1(\Omega;\R^3)$
\begin{align}  \label{eq:continuity}
& \left|\int_\Omega \TT \u: \D \v \, d\x\right|=\left|2\mu \int_\Omega \D \u:\D \v \, d\x+\lambda \int_\Omega (\dive \u) (\dive \v) \, d\x\right| \le (\lambda+2\mu) \, \|\u \, \|_{H^1} \, \|\v \, \|_{H^1} \, .
\end{align}
Estimate \eqref{eq:continuity} proves the continuity of the bilinear form
\begin{equation*}
a(\u,\v)=\int_\Omega \TT \u: \D \v \, d\x  \, , \qquad  \u,\v \in H^1(\Omega;\R^3) \, .
\end{equation*}
By Korn inequality we also see that $a(\cdot,\cdot)$ is \textit{weakly coercive}; indeed for any $0<\varepsilon<4\mu$ we have
\begin{align} \label{eq:coercive}
& a(\u,\u)+\varepsilon \|\u \, \|_{L^2}^2 \ge 2\mu \left(1-\frac{\varepsilon}{4\mu}\right) \int_\Omega |\D \u|^2 d\x
+\frac{\varepsilon}{2} \left(\frac 1C \int_\Omega |\nabla \u|^2 d\x-\int_\Omega |\u|^2 d\x\right)
+\varepsilon \int_\Omega |\u|^2 d\x  \\[5pt]
\notag   & \ge \frac{\varepsilon}{2C} \int_\Omega |\nabla \u|^2 d\x+\frac{\varepsilon}{2}  \int_\Omega |\u|^2 d\x
\ge \min\left\{\frac{\varepsilon}{2C},\frac{\varepsilon}{2} \right\} \|\u \, \|_{H^1}^2 \, .
\end{align}
On the other hand, it is easy to check that the linear functional $\Lambda:H^1(\Omega;\R^3)\to \R$ defined by
$$
\Lambda(\v)=\int_\Omega \f \cdot \v \, d\x+\int_{\partial \Omega} \g \cdot \v \, d\x
\, , \qquad \v\in H^1(\Omega;\R^3)
$$
is continuous thanks to the H\"older inequality and the classical trace inequality for $H^1$-functions. Hence, we may write
$\Lambda \in (H^{1}(\Omega;\R^3))'$.

With the notations introduced in this proof, the variational problem \eqref{eq:weak} may be written in the form
\begin{equation*}
a(\u,\v)= \langle \Lambda,\v \rangle \qquad \text{for any } \v\in H^1(\Omega;\R^3) \, .
\end{equation*}
Introducing the linear continuous operator $L:H^1(\Omega;\R^3)\to (H^1(\Omega;\R^3))'$ defined by
$$
\langle L\u,\v\rangle=a(\u,\v)  \qquad \text{for any } \u, \v \in H^1(\Omega;\R^3) \, ,
$$
we may write \eqref{eq:weak} in the form
\begin{equation} \label{eq:Lu-L}
L\u=\Lambda
\end{equation}
as an identity between elements of the dual space $(H^1(\Omega;\R^3))'$.

The next step is to introduce the following operator $R_\e:(H^1(\Omega;\R^3))'\to H^1(\Omega;\R^3)$ which maps each element
$\h \in (H^1(\Omega;\R^3))'$ into the unique solution $\w$ of the variational problem
\begin{equation*}
a(\w,\v)+\e (\w,\v)_{L^2}=\langle \h, \v\rangle  \qquad \text{for any } \v \in H^1(\Omega;\R^3) \, .
\end{equation*}
This problem admits a unique solution by the continuity and coercivity estimates \eqref{eq:continuity}, \eqref{eq:coercive} combined with the Lax-Milgram Theorem. In particular $R_\e$ is well defined and continuous. Moreover, $R_\e$ is invertible and by the Open Mapping Theorem its inverse is also continuous.

In the rest of the proof we denote by $J:H^1(\Omega;\R^3)\to (H^1(\Omega;\R^3))'$ the linear operator defined by
\begin{equation*}
\langle J\u, \v\rangle=\int_\Omega \u \cdot \v \, d\x  \qquad \text{for any } \u, \v \in H^1(\Omega;\R^3)  \, ,
\end{equation*}
which is compact as a consequence of the compact embedding $H^1(\Omega;\R^3)\subset L^2(\Omega;\R^3)$.

We now introduce on $H^1(\Omega;\R^3)$ the following scalar product
\begin{equation*}
(\u,\v)_\eps=a(\u,\v)+\e (\u,\v)_{L^2} \qquad \text{for any } \u,\v \in H^1(\Omega;\R^3) \, ,
\end{equation*}
which is equivalent to the natural scalar product of $H^1(\Omega;\R^3)$ thanks to \eqref{eq:continuity} and \eqref{eq:coercive}.

In this way we may now define the compact self-adjoint linear operator $T_\e:H^1(\Omega;\R^3)\to H^1(\Omega;\R^3)$ given by
$R_\e\circ J$ where by self-adjoint we mean
$(T_\e \u,\v)_\e=(\u,T_\e \v)_\e$  for any $\u,\v \in H^1(\Omega;\R^3)$.
Indeed, from the definition of $R_\e$, $J$ and $(\cdot,\cdot)_\e$ we see that
\begin{equation*}
(T_\e \u,\v)_\e=(\u,\v)_{L^2}=(\u,T_\e \v)_\e \qquad \text{for any }  \u,\v \in H^1(\Omega;\R^3) \, .
\end{equation*}
By definition of $R_\e$ and $J$ we have that $L=R_\e^{-1}-\e J$.
In particular $\u \in H^1(\Omega;\R^3)$ is a solution of \eqref{eq:Lu-L} if and only if \
$-\e^{-1} \, R_\e (R_\e^{-1}-\e J)\u=-\e^{-1} \, R_\e \Lambda$ \
or equivalently $T_\e \u-\e^{-1} \, \u=\w$ once we put $\w=-\e^{-1} \, R_\e \Lambda$.
Then, applying the Fredholm alternative to the operator $T_\e$ we deduce that \eqref{eq:Lu-L}, or equivalently \eqref{eq:weak}, admits a solution $\u \in H^1(\Omega;\R^3)$ if and only if
\begin{equation} \label{eq:comp-2}
\w \in \left({\rm Ker}\left(T_\eps^*-\e^{-1} \, I_{H^1}\right)\right)^\perp=\left({\rm Ker}\left(T_\eps-\e^{-1}\, I_{H^1}\right)\right)^\perp
\end{equation}
where $T_\eps^*$ denotes the adjoint operator of $T_\eps$, $I_{H^1}$ denotes the identity map in $H^1(\Omega;\R^3)$ and the orthogonal spaces are defined in the sense of the scalar product $(\cdot, \cdot)_\e$.

We observe that a function $\v\in {\rm Ker}\left(T_\eps-\e^{-1} \, I_{H^1}\right)$  if and only if
$R_\eps J \v=\e^{-1} \, \v$
and by the definition of $R_\e$ this is equivalent to
\begin{equation*}
a\left(\e^{-1}\, \v,{\bm \phi}\right)+\e \left(\e^{-1} \, \v,{\bm \phi}\right)_{L^2}=\langle J\v,{\bm \phi}\rangle
\qquad \text{for any } {\bm \phi}\in H^1(\Omega;\R^3)
\end{equation*}
and, in turn, recalling the definition of $J$ this is equivalent to $a(\v,{\bm \phi})=0$ for any ${\bm \phi}\in H^1(\Omega;\R^3)$. This shows that ${\rm Ker}\left(T_\eps-\e^{-1} \, I_{H^1}\right)=V_0$ as we deduce by \eqref{eq:V0}.

Let us proceed by proving (i)-(iv).

The proof of (i) is complete once we show that \eqref{eq:comp-2} is equivalent to condition \eqref{eq:compatibility}. Condition \eqref{eq:comp-2} is equivalent to
\begin{equation} \label{eq:final-1}
a(\w,\v)+\e (\w,\v)_{L^2}=0 \qquad \text{for any } \v\in V_0,
\end{equation}
being ${\rm Ker}\left(T_\eps-\e^{-1} \, I_{H^1}\right)=V_0$.
But $\w=-\e^{-1} \, R_\e \Lambda$ so that by definition of $R_\e$, we infer
\begin{align} \label{eq:final-2}
& a(\w,\v)+\e (\w,\v)_{L^2}=-\e^{-1} \left[a(R_\e \Lambda,\v)+\e (R_\e \Lambda,\v)_{L^2}\right]=-\e^{-1} \, \langle \Lambda,\v \rangle \qquad \text{for any } v\in H^1(\Omega;\R^3) \, .
\end{align}

Combining \eqref{eq:final-1} and \eqref{eq:final-2} we finally obtain
$\langle \Lambda,\v \rangle=0$ for any $v\in V_0$,
which is exactly \eqref{eq:compatibility} in view of the definition of the functional $\Lambda$.

For the proof of (ii) we observe that by \eqref{Hooke}, \eqref{eq:weak}, \eqref{eq:V0} and \eqref{eq:var-eig} we have for any $\v\in H^1(\Omega;\R^3)$
\begin{align*}
& \int_\Omega \TT (\u+\v_0):\D\v \, d\x=\int_\Omega \TT \u:\D\v \, d\x+\int_\Omega \TT \v_0:\D\v \, d\x =\int_\Omega \TT \u:\D\v \, d\x=\int_\Omega \f \cdot \v \, d\x+\int_{\partial\Omega} \g \cdot \v \, dS
\end{align*}
which shows that $\u+\v_0$ is a solution of \eqref{eq:weak}.

For the proof of (iii) we consider two solutions $\u$ and $\w$ of \eqref{eq:weak} and let $\v_0=\w-\u$.	
By \eqref{Hooke} and \eqref{eq:weak} we obtain
\begin{align*}
& \int_\Omega \TT \v_0 : \D \v \, d\x=\int_\Omega \TT \w : \D \v \, d\x-\int_\Omega \TT \u : \D \v \, d\x \\[5pt]
& \qquad =\int_\Omega \f \cdot \v \, d\x+\int_{\partial\Omega} \g \cdot \v \, dS-\int_\Omega \f \cdot \v \, d\x-\int_{\partial\Omega} \g \cdot \v \, dS=0  \qquad \text{for any } \v \in H^1(\Omega;\R^3)
\end{align*}
which immediately gives $\v_0\in V_0$ thanks to \eqref{eq:V0}.

Finally, let us proceed with the proof of (iv). First we prove the existence of a solution of \eqref{eq:weak} in $V_0^\perp$.

Let $\u$ be a generic solution of \eqref{eq:weak} and consider its orthogonal decomposition $\u=\u_0+\u_1 \in V_0 \oplus
V_0^\perp$ with respect to the scalar product \eqref{eq:scal-prod-T}. Then, $\u_1=\u-\u_0\in V_0^\perp$ and by part (ii) we deduce that $\u_1$ is still a solution of \eqref{eq:weak}.

Once we have proved existence, let us prove uniqueness. Let $\u, \w\in V_0^\perp$ be two solutions of  \eqref{eq:weak}. Then, on one hand we have that $\u-\w\in V_0^\perp$ and on the other hand $\u-\v\in V_0$ thanks to part (iii). Therefore, $\u-\w \in V_0\cap V_0^\perp=\{0\}$ and this readily implies $\u=\w$ thus completing the proof of (iv).
\endproof

\subsection{Proof of Proposition \ref{prop0}}\label{proofprop0}

Concerning part (i) of the Proposition we only give the proof of \eqref{eq:symmetry-1} since the proof of \eqref{eq:symmetry-2}-\eqref{eq:symmetry-3} can be obtained with a similar procedure. For any function
	$\v \in H^1(\Omega;\R^3)$ we denote by $\bar \v=(\bar v_1,\bar v_2, \bar v_2)\in H^1(\Omega;\R^3)$ the function defined by
	\begin{equation} \label{eq:reflection}
	  \bar v_1(x,y,z)=v_1(x,y,-z) \, , \ \ \
	\bar v_2(x,y,z)=v_2(x,y,-z) \, , \ \ \
	\bar v_3(x,y,z)=-v_3(x,y,-z) \, , \quad \text{$\forall \, (x,y,z)\in \Omega$.}
	\end{equation}
	Let $\u$ be the unique solution of \eqref{eq:weak-2} in $V_0^\perp$ and let
	$\bar \u$ be the corresponding function defined by \eqref{eq:reflection}.

  We start by showing that $\bar \u$ solves problem \eqref{eq:weak-2}. In doing this we show that it solves the variational problem \eqref{eq:weak} where in the present case $\f=\0 $ and $\g$ is the function defined in \eqref{eq:def-g}.
	
	By direct computation one can see that for any test function $\v \in H^1(\Omega;\R^3)$ we have
   for any $(x,y,z)\in \Omega$
	\begin{equation} \label{eq:Dw}
	(\D \bar \u: \D \v)_{|(x,y,z)}=(\D \u: \D \bar \v)_{|(x,y,-z)}  \, , \qquad
		[(\dive \, \bar \u)(\dive \, \v]_{|(x,y,z)}= [(\dive \, \u)(\dive \, \bar \v]_{|(x,y,-z)}  \, .
	\end{equation}
	
	By \eqref{eq:weak}, \eqref{Hooke}, \eqref{eq:Dw}, \eqref{eq:def-g} and a change of variables, we obtain
	\begin{align} \label{eq:u-bar}
	& 2\mu \int_{\Omega} \D \bar \u: \D \v \, d\x+\lambda \int_\Omega (\dive \, \bar \u)(\dive \, \v) \, d\x \\[5pt]
	\notag  &  \qquad =2\mu \int_{\Omega} \D \u: \D \bar \v \, d\x+\lambda \int_\Omega (\dive \, \u)(\dive \, \bar \v) \, d\x =\int_{\partial \Omega} \g \cdot \bar \v \, dS=\int_{\partial \Omega} \g \cdot \v \, dS \, .
	\end{align}
	By \eqref{eq:u-bar} we deduce that $\bar \u$ is a solution of \eqref{eq:weak} and hence a weak solution of \eqref{eq:weak-2}. We now prove that $\bar \u\in V_0^\perp$. Indeed, proceeding as in \eqref{eq:u-bar} one can easily show that
$(\bar \u,\v)_\TT=(\u,\bar \v)_\TT=0$ for any $\v \in V_0$
since $\u \in V_0^\perp$ and $\bar \v \in V_0$ whenever $\v\in V_0$, as one can deduce by \eqref{eq:van-cond}.
This completes the proof of \eqref{eq:symmetry-1}.
	
	Let us proceed with the proof of part (ii) and (iii) of the proposition.
For any $\theta\in (-2\pi,2\pi)$ we denote by $R_\theta:\R^2\to \R^2$ the anticlockwise rotation of an angle $\theta$ and by 	$A_\theta$ the associate matrix. Clearly we have that the inverse map of $R_\theta$ is given by $R_{-\theta}$ and $A_\theta^{-1}=A_{-\theta}$.
	
	We use the notation $\u=(u',u_3)\in \R^2\times \R$ with $u'=(u_1,u_2)$ and we denote by
	$$
	\nabla'u'=
	\begin{pmatrix}
	\frac{\partial u_1}{\partial x} & \frac{\partial u_1}{\partial y} \\[5pt]
	\frac{\partial u_2}{\partial x} & \frac{\partial u_2}{\partial y}
	\end{pmatrix}
	$$
	its Jacobian matrix in the $x$ and $y$ variables, and
	by $\D'u'$ the corresponding symmetric gradient given by $\frac 12 \Big(\nabla' u'+(\nabla' u')^T\Big)$; more in general, throughout this proof we will use the symbol $\nabla'$ for denoting the gradient with respect to the $x$ and $y$ variables.
	
	We now define
	\begin{equation*}
	\u_\theta(x,y,z)=\Big(R_{-\theta}\Big(u_1(R_\theta(x,y),z),u_2(R_\theta(x,y),z)\Big),u_3(R_\theta(x,y),z)\Big)
	\quad \text{for any } (x,y,z)\in \Omega \, .
	\end{equation*}
	Then, the Jacobian matrix $\nabla \u_\theta\in \R^{3\times 3}$ and in turn the matrix $\D \u_\theta$ admit a representation in terms of four blocks of dimensions $2\times 2$, $2\times 1$, $1\times 2$, $1\times 1$ respectively. We proceed directly with the representation of $\D u_\theta$:
	\begin{equation} \label{eq:blocks}
	\D \u_\theta \!=\!
	\begin{pmatrix}
	A_{-\theta} \ \D'u'\Big(R_\theta(x,y),z\Big) A_\theta &  A_{-\theta}  \frac{\partial u'}{\partial z}\Big(R_\theta(x,y),z\Big)\!+\!\left[\nabla' u_3\Big(R_\theta(x,y),z\Big) A_\theta\right]^T \\[15pt]
	\left[ A_{-\theta}  \frac{\partial u'}{\partial z}\Big(R_\theta(x,y),z\Big)\right]^T \!\!\!+\!\!\nabla' u_3\Big(R_\theta(x,y),z\Big) A_\theta            & \quad  \frac{\partial u_3}{\partial z}\Big(R_\theta(x,y),z\Big)
	\end{pmatrix}
	\end{equation}
	In the same way, for any test function $\v\in H^1(\Omega;\R^3)$ and any $\theta \in (-2\pi,2\pi)$ we may define the corresponding function $\v_\theta$. Looking at $\v$ as $(\v_{-\theta})_\theta$ and applying \eqref{eq:blocks} to $\v_{-\theta}$ we claim that for any $(x,y,z)\in \Omega$
	\begin{equation} \label{eq:sym-Du}
	\D \u_\theta(x,y,z):\D \v(x,y,z)=\D \u\Big(R_\theta(x,y),z\Big):\D \v_{-\theta}\Big(R_{\theta}(x,y),z\Big) \, .
	\end{equation}
This is a consequence of the fact that $A_\theta$ is orthogonal and the linear map $\mathcal L_\theta:\R^{2\times 2} \to \R^{2\times 2}$, $\mathcal L_\theta(X)=A_{-\theta} X A_\theta$ is an isometry in $\R^{2\times 2}$ as one can see by verifying the orthogonality of the associated matrix $M_\theta\in \R^{4\times 4}$. This implies
	\begin{align*}
	& (A_{-\theta}\, X A_\theta):Y=\mathcal L_\theta(X):Y=\mathcal L_\theta(X):\mathcal L_\theta(\mathcal L_\theta^{-1}(Y))
	=X:\mathcal L_\theta^{-1}(Y)=X:(A_\theta Y A_{-\theta})
	\end{align*}
	for any $X,Y\in \R^{2 \times 2}$. This arguments allow to treat the scalar products between the $2\times 2$ block appearing in the representation \eqref{eq:blocks}. Even easier is to treat the scalar products between the $2\times 1$ and $1\times 2$ blocks thanks to the orthogonality of $A_\theta$. This proves the claim \eqref{eq:sym-Du}.

	The invariance of the trace of a matrix $X$ under maps of the form $X\mapsto A^{-1} X A$ combined with \eqref{eq:blocks} shows that $\dive \, \u_\theta(x,y,z)=\dive \, \u\Big(R_\theta(x,y),z\Big)$ and in particular for any $(x,y,z)\in \Omega$ we have
	\begin{equation} \label{eq:sym-divu}
	(\dive \, \u_\theta(x,y,z))(\dive \, \v(x,y,z))=\left(\dive \, \u\Big(R_\theta(x,y),z\Big)\right)
	\left(\dive \, \v_{-\theta}\Big(R_{\theta}(x,y),z\Big)\right) \, .
	\end{equation}
	By \eqref{eq:weak}, \eqref{eq:def-g}, \eqref{load0}, \eqref{eq:sym-Du}, \eqref{eq:sym-divu}, two changes of variables and the definitions of $\v_{-\theta}$ and $\g$, we obtain
	\begin{align} \label{eq:inv-rot}
	& 2\mu \int_\Omega \D \u_\theta : \D \v \, d\x+\lambda \int_\Omega (\dive \, \u_\theta) (\dive \, \v) \, d\x \\[5pt]
	\notag & \quad =2\mu \int_\Omega \D \u: \D \v_{-\theta} \, d\x
	+\lambda \int_\Omega \left(\dive \, \u\right)
	\left(\dive \, \v_{-\theta}\right) \, d\x=\int_{\partial \Omega} \g \cdot \v_{-\theta} \, dS
=\int_{\partial \Omega} \g\cdot \v \, dS \, .
	\end{align}
	
	We have just proved that $\u_\theta$ is still a weak solution of \eqref{eq:weak-2}. We now show that $\u_\theta \in V_0^\perp$ as a consequence of the fact that $\u \in V_0^\perp$.
Proceeding as in \eqref{eq:inv-rot}, we infer
\begin{equation} \label{eq:scal-u-theta}
   (\u_\theta,\v)_\TT=(\u,\v_{-\theta})_\TT \qquad \text{for any } \v \in H^1(\Omega;\R^3)  \, .
\end{equation}

We need to prove that if $\v \in V_0$ then $\v_{-\theta}\in V_0$. For any $\theta \in (-2\pi,2\pi)$, let $B_\theta$ be the $3 \times 3$ matrix corresponding to an anticlockwise rotation of an angle $\theta$ around the $z$ axis. Clearly
$B_\theta$ is orthogonal and $B_\theta^{-1}=B_{-\theta}$.
With this notation we may write
\begin{equation} \label{eq:v-theta-1}
  \v_{-\theta}(\x)=B_\theta \, \v(B_{-\theta}\, \x)    \qquad \text{for any } \x\in \R^3
\end{equation}
where both $\x$ and $\v$ have to be considered vector columns in the right hand side of the identity.

If $\v \in V_0$, then by \eqref{eq:van-cond} we have that $\v$ admits the following matrix representation
\begin{equation} \label{eq:v-V0-1}
   \v(\x)=M \x+{\bm \delta}   \qquad \text{for any } \x\in \R^3
\end{equation}
where $M$ is an antisymmetric matrix and ${\bm \delta}=(\delta_1 \ \ \delta_2 \ \ \delta_3)^T$.

Combining \eqref{eq:v-theta-1} and \eqref{eq:v-V0-1} we obtain
$\v_{-\theta}(\x)=B_\theta M B_{-\theta} \, \x+B_\theta {\bm \delta}$
where the matrix $B_\theta M B_{-\theta}$ is antisymmetric since
$$
    (B_\theta M B_{-\theta})^T=B_{-\theta}^T \,  M^T B_\theta^T=B_{-\theta}^{-1} (-M) B_\theta^{-1}
    =-B_\theta M B_{-\theta} \, .
$$
This proves that also $\v_{-\theta}\in V_0$ since it admits a representation like in \eqref{eq:van-cond}.

Now, if we choose $\v \in V_0$ in \eqref{eq:scal-u-theta}, we readily see that $(\u_\theta, \v)_\TT=0$ being $\u \in V_0^\perp$ and $\v_{-\theta} \in V_0$.
This proves that $\u_\theta \in V_0^\perp$.

By the uniqueness result stated in Theorem \ref{p:existence} (iv) we infer that $\u_\theta=\u$ for any $\theta \in (-2\pi,2\pi)$.

	Now the validity of (ii) and of the first part of (iii) follows immediately from the definition of $\u_\theta$.
	
	It remains to observe that the vector field $u'$ is oriented radially in the $xy$-plane. To do this, it is sufficient to combine the identity $\u=\u_\theta$ with the identity
	$u_2(x,0,z)=0$, valid for any $a<x<b$ and $z\in \left(-\frac h 2,\frac h 2\right)$, as a consequence of \eqref{eq:symmetry-3}. \endproof
	
\subsection{Proof of Lemma \ref{lemma1}}\label{prooflemma1}
	Let us introduce the sequence of intervals $I_k:=\left(-\frac h 2+kh,\frac h 2+kh\right)$, the corresponding sequence of domains $\Omega_k:=C_{a,b}\times I_k$ and the sequence of functions $\g_k:\partial\Omega_k\to \R^3$
	\begin{equation}  \label{eq:gk}
		\g_k(x,y,z):=
		\begin{cases}
			(0,0,(-1)^{k} \, \chi_p(x,y))  & \qquad \text{if } (x,y,z)\in C_{a,b} \times \left\{-\frac h 2+kh \right\} \, ,  \\[5pt]
			(0,0,(-1)^{k+1}     \, \chi_p(x,y))  & \qquad \text{if } (x,y,z)\in  C_{a,b} \times \left\{\frac h 2+kh \right\} \, ,
			\\[5pt]
			(0,0,0)                       & \qquad \text{if } (x,y,z)\in \partial C_{a,b} \times I_k  \, .
		\end{cases}
	\end{equation}
	We know that the original function $\u$ is a weak solution of problem \eqref{eq:weak-2} in the sense that
	\begin{equation} \label{eq:var-3}
		\int_{\Omega} \TT \u: \D \v \, d\x=\int_{\partial\Omega} \g \cdot \v \, dS \qquad \text{for any } \v\in H^1(\Omega;\R^3) \, .
	\end{equation}
	We need to find, starting from \eqref{eq:var-3}, the equation solved, in the sense of distributions, by the periodic extension.
	First of all, we observe that by \eqref{eq:sym-ext-123}, \eqref{eq:gk}, \eqref{eq:var-3} and some computations, we have
	\begin{equation} \label{eq:var-4}
		\int_{\Omega_k} \TT \u: \D \v \, d\x=\int_{\partial \Omega_k} \g_k \cdot \v \, dS \qquad \text{for any } \v\in H^1(\Omega_k;\R^3) \, .
	\end{equation}
	Now, letting ${\bm \phi}=(\phi_1,\phi_2,\phi_3)\in \mathcal D(C_{a,b}\times \R;\R^3)$, by \eqref{eq:var-4} we infer
	\begin{align*}
		& \int_{C_{a,b}\times \R} \TT \u : \D {\bm \phi} \, d\x=\sum_{k\in \Z} \, \int_{\Omega_k} \TT  \u : \D {\bm \phi} \, d\x
		=\sum_{k\in \Z} \, \int_{\partial \Omega_k} \g_k \cdot {\bm \phi} \, dS  \\[10pt]
		& \qquad =\sum_{k\in \Z} 2 (-1)^{k+1} \int_{C_{a,b}} \chi_p(x,y)\,  \phi_3\left(x,y,\tfrac h 2+kh \right) dxdy \, .
	\end{align*}
	
	This proves \eqref{eq:distr-equation},	where ${\bm \Lambda}$ is the distribution defined by
	\begin{equation} \label{eq:def-Lambda-1}
		\langle {\bm \Lambda}, {\bm \phi} \rangle:= \int_{C_{a,b}}  \chi_p(x,y)\,  \sum_{k\in \Z} 2(-1)^{k+1}\phi_3\left(x,y,\tfrac h 2+kh \right) dxdy
	\end{equation}
	for any ${\bm \phi}=(\phi_1,\phi_2,\phi_3)\in \mathcal D(C_{a,b}\times \R;\R^3)$.
	
	The distribution ${\bm \Lambda}$ admits a sort of factorization as a product of a function in the variables $x$ and $y$ and
	of a distribution acting on functions of the variable $z$:
	\begin{equation*}
		{\bm \Lambda}=(\Lambda_1,\Lambda_2,\Lambda_3)=\left(0,0,2\, \chi_p \,  \sum_{k\in \Z} (-1)^{k+1} \ \delta_{\frac h 2+kh}\right)
	\end{equation*}
	where $\Lambda_1, \Lambda_2, \Lambda_3 \in \mathcal D'(C_{a,b}\times \R;\R)$ are the \textit{scalar distributions} defined by
	\begin{equation*}
		\langle \Lambda_i,\phi\rangle:=\langle {\bm \Lambda},\phi\,  {\bf e}_i\rangle \qquad \text{for any } \phi \in \mathcal D(C_{a,b}\times \R;\R) \, ,
	\end{equation*}
	with ${\bf e}_1=\i$, ${\bf e}_2=\j$, ${\bf e}_3=\k$, and $\delta_{\frac h 2+kh}$ are Dirac delta distributions concentrated at $z=\frac h 2+kh$.
	
	Expanding in Fourier series the periodic distribution $\sum_{k\in \Z} \ 2(-1)^{k+1} \ \delta_{\frac h 2+kh}$ we obtain \eqref{eq:distr-expansion}, where the Fourier series converges in the sense of distributions. For more details on this convergence see the arguments introduced in subsection \ref{ss:main-teo}.\endproof

\subsection{Proof of Lemma \ref{lemma2}}\label{prooflemma2}
First of all we insert \eqref{eq:comp-expansion} into \eqref{eq:distr-equation}; recalling the Hooke's law \eqref{Hooke} and exploiting \eqref{eq:distr-expansion}, we obtain

\begin{equation}\label{pbcil2}
	\begin{cases}
		-\mu \Delta\varphi_k^1+\mu \dfrac{\pi^2 k^2}{h^2} \, \varphi_k^1-(\lambda+\mu)\bigg[\dfrac{\partial^2 \varphi_k^1}{\partial x^2}+\dfrac{\partial^2 \varphi_k^2}{\partial x \partial y}+\dfrac{\pi k}h\dfrac{\partial \varphi_k^3}{\partial x}\bigg]=0 \quad&\text{in }  C_{a,b} \, , \\[12pt]
		-\mu \Delta\varphi_k^2+\mu \dfrac{\pi^2 k^2}{h^2} \, \varphi_k^2-(\lambda+\mu) \bigg[\dfrac{\partial^2 \varphi_k^1}{\partial x\partial y}+\dfrac{\partial^2 \varphi_k^2}{\partial y^2}+\dfrac{\pi k}h\dfrac{\partial \varphi_k^3}{\partial y}\bigg]=0 \quad &\text{in } C_{a,b} \, , \\[12pt]
		-\mu \Delta\varphi_k^3+\mu \dfrac{\pi^2 k^2}{h^2} \, \varphi_k^3+\dfrac{\pi k}h (\lambda+\mu) \bigg[\dfrac{\partial \varphi_k^1}{\partial x}+\dfrac{\partial \varphi_k^2}{\partial y}+\dfrac{\pi k}h \varphi_k^3\bigg]=\Psi_k \quad&\text{in }  C_{a,b} \, ,
	\end{cases}
\end{equation}
where the forcing term is defined in \eqref{eq:Psi-k}.
We observe that in \eqref{pbcil2}, the operator $\Delta$ stands for the Laplace operator in the variables $x$ and $y$, i.e.
$\Delta=\partial^2/\partial x^2+\partial^2/\partial y^2$.

Putting $\Phi_k:=(\varphi_k^1,\varphi_k^2)$ and $\bar \n\in \R^2$ the outward unit normal to $\partial C_{a,b}$, system \eqref{pbcil2} may be rewritten in the following form

\begin{equation} \label{eq:alt-ver}
	\begin{cases}
		-\mu \Delta \Phi_k+\mu \frac{\pi^2 k^2}{h^2} \Phi_k-(\lambda+\mu) \nabla (\dive \, \Phi_k)-(\lambda+\mu) \frac{\pi k}h \nabla \varphi_k^3=\0
		& \qquad \text{in }  C_{a,b} \, , \\[8pt]
		-\mu \Delta \varphi_k^3+\mu \frac{\pi^2 k^2}{h^2} \varphi_k^3+\frac{\pi k}h(\lambda+\mu) \left(\dive \, \Phi_k+\frac{\pi k}h \varphi_k^3\right)=\Psi_k
		& \qquad \text{in }  C_{a,b} \, ,
	\end{cases}
\end{equation}
or equivalently in the following form
\begin{equation} \label{eq:2D}
	\begin{cases}
		-\dive(\lambda (\dive \, \Phi_k)I+2\mu \D \Phi_k)+\mu \frac{\pi^2 k^2}{h^2} \Phi_k-(\lambda+\mu) \frac{\pi k}h \nabla \varphi_k^3=\0
		& \qquad \text{in }  C_{a,b} \, , \\[7pt]
		-\mu \Delta \varphi_k^3+(\lambda+2\mu) \frac{\pi^2 k^2}{h^2} \varphi_k^3+(\lambda+\mu)\frac{\pi k}h \,\dive \, \Phi_k=\Psi_k
		& \qquad \text{in }  C_{a,b} \, ,
	\end{cases}
\end{equation}
where $\D$ represents here the symmetric gradient in the two-dimensional case and $I$ is the $2\times 2$ identity matrix.

We also recall that by \eqref{eq:weak-2}, $(\TT  \u)\n=\0 $ on $\partial  C_{a,b} \times \R$ so that by the Hooke's law \eqref{Hooke} we obtain
\begin{equation*}
	\begin{cases}
		(\lambda+2\mu) x \dfrac{\partial u_1}{\partial x}+\lambda x \dfrac{\partial u_2}{\partial y}
		+\lambda x \dfrac{\partial u_3}{\partial z}+\mu y \left(\dfrac{\partial u_1}{\partial y}+\dfrac{\partial u_2}{\partial x}\right)=0 & \qquad \text{on } \partial  C_{a,b}\times \R \, , \\[15pt]
		\lambda y\dfrac{\partial u_1}{\partial x}+(\lambda+2\mu) y \dfrac{\partial u_2}{\partial y}
		+\lambda y \dfrac{\partial u_3}{\partial z}+\mu x\left(\dfrac{\partial u_1}{\partial y}+\dfrac{\partial u_2}{\partial x}\right)=0 & \qquad \text{on } \partial  C_{a,b}\times \R \, , \\[15pt]
		\mu x \left(\dfrac{\partial u_1}{\partial z}+\dfrac{\partial u_3}{\partial x}\right)+\mu y \left(\dfrac{\partial u_2}{\partial z}+\dfrac{\partial u_3}{\partial y}\right)=0 & \qquad \text{on } \partial  C_{a,b}\times \R \, ,
	\end{cases}
\end{equation*}
and by \eqref{eq:comp-expansion} we obtain
\begin{equation} \label{eq:BC-2D}
	\begin{cases}
		\lambda(\dive\, \Phi_k) \bar \n+2\mu (\D\Phi_k)\bar \n+\lambda \frac{\pi k}h \varphi_k^3\,  \bar \n=\0 & \qquad \text{on } \partial  C_{a,b} \, , \\[10pt]
		\mu \nabla \varphi_k^3 \cdot \bar \n-\mu \frac{\pi k}h \, \Phi_k\cdot \bar \n=0 & \qquad \text{on } \partial  C_{a,b} \, .
	\end{cases}
\end{equation}
Let us derive the weak formulation of \eqref{eq:alt-ver}-\eqref{eq:BC-2D}.
Testing \eqref{eq:2D} with $(w^1,w^2,w^3)$, putting $W=(w^1,w^2)$ and integrating by parts we obtain
\begin{align} \label{eq:big-var}
	& -\int_{\partial  C_{a,b}} \left[\left(\lambda (\dive \, \Phi_k)I+2\mu \D \Phi_k\right)\bar \n\right]\cdot W \, ds
	+\int_{ C_{a,b}} \left(\lambda (\dive \, \Phi_k)I+2\mu \D \Phi_k\right) : \nabla W \, dxdy  \\[5pt]
	\notag & \quad + \mu \frac{\pi^2 k^2}{h^2} \int_{ C_{a,b}}\!\! \Phi_k \cdot W \, dxdy-\mu \frac{\pi k}{h} \int_{ C_{a,b}} \!\! \nabla \varphi_k^3 \cdot W \, dxdy
	-\lambda \frac{\pi k}{h} \int_{\partial  C_{a,b}} \!\! \varphi_k^3 \, \bar \n \cdot W \, ds+\lambda \frac{\pi k}{h} \int_{ C_{a,b}} \!\! \varphi_k^3 \, \dive \, W \, dxdy \\[5pt]
	\notag & \quad -\int_{\partial  C_{a,b}} \mu \frac{\partial \varphi_k^3}{\partial \bar \n} \, w^3 ds +\mu \int_{ C_{a,b}} \nabla \varphi_k^3 \cdot \nabla w^3 \, dxdy+(\lambda+2\mu) \frac{\pi^2 k^2}{h^2} \int_{ C_{a,b}} \varphi_k^3 \, w^3 \, dxdy  \\[5pt]
	\notag  & \quad +\mu \frac{\pi k}{h} \int_{\partial  C_{a,b}} w^3 \Phi_k \cdot \bar \n \, ds-\mu \frac{\pi k}{h} \int_{ C_{a,b}} \Phi_k \cdot \nabla w^3\, dxdy+ \lambda \frac{\pi k}{h} \int_{ C_{a,b}} \dive \, \Phi_k \, w^3 \, dxdy=\int_{ C_{a,b}} \Psi_k \, w^3 \, dxdy \, .
\end{align}
We observe that by \eqref{eq:BC-2D} the boundary integrals in \eqref{eq:big-var} disappear; on the other hand collecting the double integrals and recalling that $\D \Phi_k : \nabla W=\D \Phi_k : \D W$, we may write \eqref{eq:big-var} in the form
\begin{align} \label{eq:big-var-2}
	&    2\mu \int_{ C_{a,b}} \D \Phi_k : \D W \, dxdy+\lambda \int_{ C_{a,b}} \left( \dive \, \Phi_k+\tfrac{\pi k}h \varphi_k^3 \right)
	\left( \dive \, W+\tfrac{\pi k}h w^3\right) dxdy  \\[7pt]
	\notag &   
	+\mu \int_{ C_{a,b}} (\nabla \varphi_k^3-\tfrac{\pi k}h \Phi_k)\cdot (\nabla w^3-\tfrac{\pi k}h W) dxdy
	+2\mu \tfrac{\pi^2 k^2}{h^2} \int_{ C_{a,b}} \varphi_k^3w^3dxdy\!  =\!\int_{ C_{a,b}} \!\Psi_k w^3  dxdy
\end{align}
for any $\w\in H^1( C_{a,b};\R^3)$, where $W=(w^1,w^2)$. This represents
the weak form of \eqref{eq:alt-ver}-\eqref{eq:BC-2D}.

For any $k\ge 2$ even we observe that, by \eqref{eq:Psi-k}, $\varphi_k^1\equiv \varphi_k^2\equiv \varphi_k^3 \equiv 0$ in $C_{a,b}$, as one can deduce by testing \eqref{eq:big-var-2} with $(w^1,w^2,w^2)=(\varphi_k^1,\varphi_k^2,\varphi_k^3)$.

For any $k\ge 1$ odd, we define the following bilinear form
\begin{align} \label{eq:bar-a}
	& \bar a_k({\bm \varphi},\w):= 2\mu \int_{ C_{a,b}} \D \Phi : \D W \, dxdy+\lambda \int_{ C_{a,b}} \left( \dive \, \Phi+\tfrac{\pi k}{h}\,  \varphi_k^3 \right)
	\left( \dive \, W+\tfrac{\pi k}{h}\,  w^3\right) dxdy  \\[5pt]
	\notag  & \qquad +\mu \int_{ C_{a,b}} (\nabla \varphi^3-\tfrac{\pi k}{h} \, \Phi)\cdot (\nabla w^3-\tfrac{\pi k}{h} \, W) \, dxdy
	+2\mu \tfrac{\pi^2 k^2}{h^2} \int_{ C_{a,b}} \varphi^3 \, w^3 \, dxdy \quad \text{for any } {\bm \varphi}, \w \in H^1( C_{a,b};\R^3)
\end{align}
where ${\bm \varphi}=(\varphi^1,\varphi^2,\varphi^3)$, $\Phi:=(\varphi^1,\varphi^2)$, $\w=(w^1,w^2,w^3)$ and $W=(w^1,w^2)$.

For the uniqueness issue we claim that for any $\varepsilon>0$ there exists $C_\e>0$ such that
\begin{equation} \label{eq:coercive-22}
	\bar a_k({\bm \varphi},{\bm \varphi})+\eps \|{\bm \varphi}\|_{L^2}^2 \ge C_\e \|{\bm \varphi}\|_{H^1}^2
	\qquad \text{for any } {\bm \varphi} \in H^1( C_{a,b};\R^3) \, .
\end{equation}

Suppose by contradiction that there exists $\e>0$ such that for any $m\ge 1$ there exists ${\bm \varphi}_m \in H^1( C_{a,b};\R^3)$ such that
\begin{equation} \label{eq:negation}
	\bar a_k({\bm \varphi}_m,{\bm \varphi}_m)+\eps \|{\bm \varphi}_m\|_{L^2}^2 \le \frac 1m \|{\bm \varphi}_m\|_{H^1}^2
	\, .
\end{equation}

Up to normalization, it is not restrictive to assume that the sequence $\{{\bm \varphi}_m\}$ satisfies $\|{\bm \varphi}_m\|_{H^1}=1$ for any $m\ge 1$, so that by \eqref{eq:bar-a} and \eqref{eq:negation} we infer
\begin{align} \label{eq:convergences}
	& {\bm \varphi}_m\to \0 \quad \text{in } L^2( C_{a,b};\R^3) \, , \qquad \int_{ C_{a,b}} |\D \Phi_m|^2 dxdy\to 0
	\, , \qquad \nabla \varphi_m^3-k\Phi_m \to \0 \quad \text{in } L^2( C_{a,b};\R^2)
\end{align}
as $m\to +\infty$.
Applying \eqref{t:Korn} in the two-dimensional case we obtain
\begin{equation*}
	\int_{ C_{a,b}} |\nabla \Phi_m|^2 dxdy\le C \left(\int_{ C_{a,b}} |\D \Phi_m|^2 dxdy+\int_{ C_{a,b}} |\Phi_m|^2 \, dxdy\right)
\end{equation*}
for some constant $C>0$. This, combined with \eqref{eq:convergences}, proves that
\begin{equation*}
	\nabla \Phi_m \to \0 \quad \text{in } L^2( C_{a,b};\R^{2\times 2}) \, , \qquad \nabla \varphi_m^3 \to \0
	\qquad \text{in } L^2( C_{a,b};\R^2)
\end{equation*}
and, in turn, that ${\bm \varphi}_m\to \0 $ in $H^1( C_{a,b};\R^3)$. This contradicts the assumption $\|{\bm \varphi}_m\|_{H^1}=1$. We have completed the proof of the claim \eqref{eq:coercive-22}.

Thanks to \eqref{eq:coercive-22}, we may proceed as in the proof of Theorem \ref{p:existence} and apply the Fredholm alternative to show that \eqref{eq:2D}-\eqref{eq:BC-2D} admits a solution if and only if
\begin{equation} \label{eq:cond-comp}
	\int_{ C_{a,b}} \Psi_k \, w^3 \, dxdy=0  \qquad \text{for any } \w=(w^1,w^2,w^3)\in \bar V_k
\end{equation}
where $\bar V_k:=\{\w \in H^1( C_{a,b};\R^3): \bar a_k(\w,\v)=0 \ \text{for any } \v\in H^1( C_{a,b};\R^3) \}$.
Testing the variational identity in the definition of $\bar V_k$ with $\v=\w$, we readily see that for any $k\ge 1$ we $\bar V_k=\{\0 \}$ and hence, condition \eqref{eq:cond-comp} is always satisfied. This completes the proof of the lemma.  \endproof

\subsection{Proof of Lemma \ref{p:prop-ex-un-rad}} \label{ss:p-33} 

Before proceeding with the proof of the lemma, we devote the first part of this subsection to show that the functions $Y_k$ and $Z_k$ introduced in \eqref{eq:Yk-Zk} really satisfy \eqref{pbcil3}.

In order to simplify the notations we denote by $Y$ and $Z$ the unknown functions, omitting the index $k$.
	Testing \eqref{eq:big-var-2} with a test function $(w^1,w^2,w^3)$ admitting in polar coordinates the following representation
	\begin{equation*}
	w^1(\rho,\theta)=H(\rho) \cos \theta \, , \quad
	w^2(\rho,\theta)=H(\rho) \sin \theta \, , \quad
	w^3(\rho,\theta)=K(\rho) \, ,
	\end{equation*}
by \eqref{eq:Yk-Zk} we obtain
\small
	\begin{align}  \label{eq:var-rad-0}
	& 2\mu \int_{ a}^{ b} \left[\rho Y'(\rho)H'(\rho)+\frac{Y(\rho)H(\rho)}{\rho}\right] d\rho
	+\lambda \int_{ a}^{ b} \rho \left[Y'(\rho)+\frac{Y(\rho)}{\rho}+\frac{\pi k}{h} \, Z(\rho)\right] \left[H'(\rho)+\frac{H(\rho)}{\rho}+\frac{\pi k}{h} \, K(\rho)\right] d\rho \\[5pt]
	\notag  & \quad +\mu \int_{ a}^{ b} \rho \left[Z'(\rho)-\frac{\pi k}{h} \, Y(\rho)\right]\left[K'(\rho)-\frac{\pi k}{h} \, H(\rho)\right] d\rho
	+2\mu \frac{\pi^2 k^2}{h^2} \int_{ a}^{ b} \rho Z(\rho)K(\rho) \, d\rho=\int_{ a}^{ b} \rho \Psi_k(\rho) K(\rho) \, d\rho
	\end{align}
\normalsize
	with obvious meaning of the notation $\Psi_k(\rho)$ being it a radial function.

	
	Collecting in a proper way the terms of \eqref{eq:var-rad-0}, we may rewrite it in the form
	\begin{align}  \label{eq:var-rad}
	& \int_{ a}^{ b} \left[(\lambda+2\mu)\rho Y'(\rho)+\lambda Y(\rho)+\lambda \frac{\pi k}{h} \,\rho Z(\rho)\right]H'(\rho) \, d\rho \\[5pt]
	\notag & \quad    +\int_{ a}^{ b} \left[\lambda Y'(\rho)+(\lambda+2\mu)\frac{Y(\rho)}{\rho}+\mu \frac{\pi^2 k^2}{h^2} \rho Y(\rho)-\mu \frac{\pi k}{h} \, \rho Z'(\rho)+\lambda \frac{\pi k}{h} \, Z(\rho)\right] H(\rho) \, d\rho \\[5pt]
	\notag & \quad +\int_{ a}^{ b} \left[\mu \rho Z'(\rho)-\mu \frac{\pi k}{h} \, \rho Y(\rho)\right]K'(\rho) \, d\rho
	+\int_{ a}^{ b} \left[(\lambda+2\mu) \frac{\pi^2 k^2}{h^2} \rho Z(\rho)+\lambda \frac{\pi k}{h} \, \rho Y'(\rho)+\lambda \frac{\pi k}{h} \, Y(\rho)\right] K(\rho) \, d\rho \\[5pt]
	& \quad =\int_{ a}^{ b} \rho \Psi_k(\rho) K(\rho) \, d\rho \notag
	\end{align}
	
	Integrating by parts the terms in \eqref{eq:var-rad} containing $H'(\rho)$ and $K'(\rho)$, we see that \eqref{eq:var-rad} is the variational formulation of \eqref{pbcil3}.

Let us proceed now with the proof of the lemma which is the main point of this section. Actually, we give here only a sketch of the proof since it essentially follows the ideas already introduced in the proof of Lemma \ref{lemma2}.

	About the uniqueness issue, on the space $H^1( a,  b;\R^2)$ it sufficient to define the bilinear form
	$$
	b_k:H^1( a, b;\R^2)\times H^1( a, b;\R^2)\to \R
	$$
	corresponding to the left hand side of \eqref{eq:var-rad-0} and prove for it an estimate of the type \eqref{eq:coercive-22}.
	
	Then, following again the proof of Lemma \ref{lemma2}, one finds that the compatibility condition for $\Psi_k$ is given by
	\begin{equation} \label{eq:compat-rad}
	\int_{ a}^{ b} \rho \Psi_k(\rho) K(\rho) \, d\rho=0
	\end{equation}
	for any $(H,K)\in H^1( a, b;\R^2)$ satisfying $b_k\Big((H,K),(H,K)\Big)=0$. A simple check shows that $(H,K)\equiv (0,0)$ so that \eqref{eq:compat-rad} is trivially satisfied.

The Fredholm alternative then implies the existence of a solution.
\endproof

\subsection{Proof of Lemma \ref{lemma3}}\label{prooflemma3}

We omit for simplicity the dependence from the index $k$ in the unknowns $Y_k$ and $Z_k$.
For more clarity we divide the construction of this representation of $Y$ and $Z$ into different steps each of them is contained in the next subsections.

\subsubsection{The solution of the homogeneous system}
We consider the homogeneous version of the system in \eqref{pbcil3}

\begin{equation}\label{eq:hom}
\begin{cases}
Y''(\rho)+\dfrac{Y'(\rho)}{\rho}-\dfrac{Y(\rho)}{\rho^2}-\alpha \, k^2 \, Y(\rho)+\beta \, k Z'(\rho)=0 \quad & \rho>0 \, , \\[12pt]
Z''(\rho)+\dfrac{Z'(\rho)}{\rho}-\gamma \, k^2 Z(\rho)-\delta \, k\bigg[Y'(\rho)+\dfrac{Y(\rho)}{\rho}\bigg]=0  \quad & \rho>0 \, ,
\end{cases}
\end{equation}
where we put for simplicity
\begin{equation*}
\alpha=\frac{\pi^2 \mu}{h^2(\lambda+2\mu)} \, , \quad \beta=\frac{\pi(\lambda+\mu)}{h(\lambda+2\mu)} \, ,
\quad \gamma=\frac{\pi^2(\lambda+2\mu)}{h^2 \mu} \, , \quad \delta=\frac{\pi(\lambda+\mu)}{h \mu} \, .
\end{equation*}
We look for a solution admitting the following expansion
\begin{equation}\label{expansion}
\begin{cases}
{\ds Y(\rho)=\sum_{n=-1}^{+\infty} a_n \, \rho^n+(\ln \rho) \sum_{n=0}^{+\infty} b_n \, \rho^n} \, , \\[15pt]
{\ds Z(\rho)=\sum_{n=0}^{+\infty} c_n \, \rho^n+(\ln \rho) \sum_{n=0}^{+\infty} d_n \, \rho^n} \, .
\end{cases}
\end{equation}

Inserting the representation \eqref{expansion} in the system \eqref{eq:hom}, we obtain for each of the two equations the following identities:

\begin{equation} \label{eq:expansion-sol-1}
\begin{tabular}{l}
${\ds \sum_{n=-1}^{+\infty} n(n-1) a_n \, \rho^n+\sum_{n=0}^{+\infty} (n-1) b_n \, \rho^n
	+\sum_{n=0}^{+\infty} n b_n \, \rho^n+(\ln \rho) \sum_{n=0}^{+\infty} n(n-1) b_n \, \rho^n}$ \\[15pt]
${\ds \qquad +\sum_{n=-1}^{+\infty} n a_n \, \rho^n+\sum_{n=0}^{+\infty} b_n \, \rho^n
	+(\ln \rho) \sum_{n=0}^{+\infty} n b_n \, \rho^n }$ \\[15pt]
${\ds \qquad -\alpha k^2 \sum_{n=1}^{+\infty} a_{n-2} \, \rho^n
	-\alpha k^2 (\ln\rho) \sum_{n=2}^{+\infty} b_{n-2} \, \rho^n-\sum_{n=-1}^{+\infty} a_{n} \, \rho^n
	-(\ln\rho) \sum_{n=0}^{+\infty} b_{n} \, \rho^n }$  \\[15pt]
${\ds \qquad +\beta k \sum_{n=1}^{+\infty} (n-1) c_{n-1} \, \rho^n+\beta k \sum_{n=1}^{+\infty} d_{n-1} \, \rho^n
	+\beta k(\ln \rho)\sum_{n=1}^{+\infty} (n-1) d_{n-1} \, \rho^n=0} \, ,$ \\[20pt]
\end{tabular}
\end{equation}

\begin{equation}  \label{eq:expansion-sol-2}
\begin{tabular}{l}
${\ds \sum_{n=0}^{+\infty} n(n-1) c_n \, \rho^n+\sum_{n=0}^{+\infty} (n-1) d_n \, \rho^n
	+\sum_{n=0}^{+\infty} n d_n \, \rho^n+(\ln \rho) \sum_{n=0}^{+\infty} n(n-1) d_n \, \rho^n}$ \\[15pt]
${\ds \qquad +\sum_{n=0}^{+\infty} n c_n \, \rho^n+\sum_{n=0}^{+\infty} d_n \, \rho^n
	+(\ln \rho) \sum_{n=0}^{+\infty} n d_n \, \rho^n-\gamma k^2 \sum_{n=2}^{+\infty} c_{n-2} \, \rho^n
	-\gamma k^2 (\ln\rho) \sum_{n=2}^{+\infty} d_{n-2} \, \rho^n }$ \\[15pt]
${\ds \qquad -\delta k \sum_{n=0}^{+\infty} (n-1) a_{n-1} \, \rho^n-\delta k \sum_{n=1}^{+\infty} b_{n-1} \, \rho^n
	-\delta k (\ln \rho)\sum_{n=1}^{+\infty} (n-1) b_{n-1} \, \rho^n}$ \\[15pt]
${\ds -\delta k \sum_{n=0}^{+\infty} a_{n-1} \, \rho^n-\delta k (\ln \rho) \sum_{n=1}^{+\infty} b_{n-1} \, \rho^n=0} \, . $
\end{tabular}
\end{equation}
To determine the values of the coefficients $a_n,b_n,c_n,d_n$ we need an iterative scheme starting from the values of the coefficients $a_{-1}, a_0, a_1, b_0, b_1, c_0, c_1, d_0, d_1$.
The values of these nine parameters have to be determined collecting the coefficients of the terms $\rho^{-1}, \rho^0, \rho^0 \ln \rho, \rho, \rho\ln \rho$ appearing in \eqref{eq:expansion-sol-1}-\eqref{eq:expansion-sol-2}  and equating them to zero.

As a result of this procedure we obtain the following constraint:
\begin{equation} \label{eq:varie}
\begin{cases}
a_0=b_0=c_1=d_1=0 \, ,  \\
2b_1+\beta kd_0=\alpha k^2 a_{-1} \, .
\end{cases}
\end{equation}
Among the left five parameters $a_{-1}, a_1, b_1, c_0, d_0$ that may be possibly different from zero, $a_1, c_0$ and two among $a_{-1}, b_1, d_0$ can be chosen arbitrarily, while the remaining one is determined by the equation in the second line of \eqref{eq:varie}; for example, we may choose arbitrarily $a_{-1}, a_1, b_1, c_0$ and put $d_0=\frac{\alpha k}{\beta} \, a_{-1}-\frac{ 2}{ \beta k}\, b_1$.

In particular, we are interested in finding the general solution of \eqref{eq:hom} as a linear combination of four linearly independent special solutions, denoted by ${\bm \Upsilon}^j=(Y^j,Z^j)$ with $j=1,\dots,4$.
A possible choice for the independent solutions is given respectively by the assumption on the following combinations of coefficients:
\begin{equation}\label{eq:cases}
\begin{split}
&{\bm \Upsilon}^1:\,\,(a_{-1},a_1,b_1,c_0)=(1,0,0,0) \, , \qquad {\bm \Upsilon}^2:\,\,(a_{-1},a_1,b_1,c_0)=(0,1,0,0) \, ,\\
&{\bm \Upsilon}^3:\,\,(a_{-1},a_1,b_1,c_0)=(0,0,1,0) \, , \qquad {\bm \Upsilon}^4:\,\,(a_{-1},a_1,b_1,c_0)=(0,0,0,1).
\end{split}
\end{equation}

By \eqref{eq:expansion-sol-1}-\eqref{eq:expansion-sol-2} we deduce the following linear system in the unknowns $a_n, b_n, c_{n-1}, d_{n-1}$ with data expressed in terms of $a_{n-2}, b_{n-2}, c_{n-3}, d_{n-3}$:
\begin{equation} \label{eq:linear-system-0}
\begin{cases}
(n^2-1) a_n+2n b_n+\beta k(n-1) c_{n-1}+\beta k d_{n-1}=\alpha k^2 a_{n-2} \\
(n^2-1) b_n +\beta k(n-1) d_{n-1}=\alpha k^2 b_{n-2} \\ 
(n-1)^2 c_{n-1}+2(n-1) d_{n-1}=\delta k(n-1) a_{n-2}+\delta k b_{n-2}+\gamma k^2 c_{n-3} \\ 
(n-1)^2 d_{n-1}=\delta k (n-1) b_{n-2}+\gamma k^2 d_{n-3} 
\end{cases}\qquad (n\geq 3).
\end{equation}

We observe that the matrix of coefficients associated to system \eqref{eq:linear-system-0} is given by

\begin{equation*}
\begin{pmatrix}
n^2-1 & 2n & \beta(n-1)k & \beta k      \\[4pt]
0&n^2-1 & 0           & \beta(n-1)k \\[4pt]
0     & 0  & (n-1)^2    & 2(n-1)       \\[4pt]
0     & 0  & 0          & (n-1)^2
\end{pmatrix}
\end{equation*}
whose determinant is given by $(n-1)^6 (n+1)^2\neq 0$, thus showing that the system is not singular for $n\ge 2$ and hence admits a unique solution.

With the restriction $n\ge 3$ the coefficients $a_2$, $b_2$ remained excluded, but their calculation can be obtained from the first two equations of \eqref{eq:linear-system-0} by choosing $n=2$;
this gives $a_2=b_2=0$.

The linear independence of ${\bm \Upsilon}^1,  {\bm \Upsilon}^2,  {\bm \Upsilon}^3,  {\bm \Upsilon}^4$ can be verified by looking at the asymptotic behavior of $Y^j(\rho)$, $j=1,2,3,4$ as $\rho\to 0^+$ in the four cases \eqref{eq:cases}:
\begin{align*}
& \text{{\bf case 1:}} \quad Y^1(\rho)\sim \rho^{-1} \quad \text{as $\rho\to 0^+$}\, ; \qquad \text{{\bf case 2:}} \quad Y^2(\rho)\sim \rho \quad \text{as $\rho\to 0^+$}\, ;\\[4pt]
& \text{{\bf case 3:}} \quad Y^3(\rho)\sim \rho \ln \rho \quad \text{as $\rho\to 0^+$}\, ;  \qquad \text{{\bf case 4:}} \quad Y^4(\rho)=O(\rho^2 \ln \rho) \quad \text{as $\rho\to 0^+$} \, .
\end{align*}

\begin{remark}
	We observe that, after a suitable scaling, the dependence of system \eqref{eq:hom} from the parameter $k$ can be dropped:
	given a solution $(Y,Z)$ of \eqref{eq:hom}, we may define the functions $\widetilde Y(t)=Y\left(\frac{h}{\pi k} \, t\right)$ and $\widetilde Z(t)=Z\left(\frac{h}{\pi k} \, t\right)$ in such a way that the couple $(\widetilde Y,\widetilde Z)$ solves system
	
	\begin{equation}\label{eq:hom-scaled}
	\begin{cases}
	\widetilde Y''(t)+\dfrac{\widetilde Y'(t)}{t}-\dfrac{\widetilde Y(t)}{t^2}-\widetilde\alpha  \, \widetilde Y(t)+\widetilde\beta \, \widetilde Z'(t)=0 \quad & t>0 \, , \\[12pt]
	\widetilde Z''(t)+\dfrac{\widetilde Z'(t)}{t}-\widetilde\gamma \, \widetilde Z(t)-\widetilde\delta \, \bigg[\widetilde Y'(t)+\dfrac{\widetilde Y(t)}{t}\bigg]=0  \quad & t>0 \, ,
	\end{cases}
	\end{equation}
	where $\widetilde \alpha=\mu/(\lambda+2\mu)$, $\widetilde \beta=(\lambda+\mu)/(\lambda+2\mu)$, $\widetilde \gamma=(\lambda+2\mu)/\mu$ and $\widetilde \delta=(\lambda+\mu)/\mu$.
\end{remark}

\subsubsection{The particular solution}

We write the nonhomogeneous system in the matrix form

\begin{equation} \label{eq:matrix-form}
\begin{pmatrix}
Y(\rho) \\[4pt]
Y'(\rho) \\[4pt]
Z(\rho) \\[4pt]
Z'(\rho)
\end{pmatrix}'
=\begin{pmatrix}
0 & 1 & 0 & 0 \\[4pt]
\frac{1}{\rho^2}+\alpha k^2 & -\frac 1\rho & 0 & -\beta k \\[4pt]
0 & 0 & 0 & 1 \\[4pt]
\frac{\delta k}{\rho} & \delta k & \gamma k^2 & -\frac 1\rho
\end{pmatrix}
\begin{pmatrix}
Y(\rho) \\[4pt]
Y'(\rho) \\[4pt]
Z(\rho) \\[4pt]
Z'(\rho)
\end{pmatrix}
+
\begin{pmatrix}
0 \\[4pt]
0 \\[4pt]
0 \\[4pt]
-\frac 1\mu \Psi_k(\rho)
\end{pmatrix}  \, , \qquad \rho>0 \, ,
\end{equation}
where the function $\Psi_k=\Psi_k(\rho)$ is extended trivially outside the interval $( a,  b)$.

Maintaining the order of the components, we may write the Wronskian matrix associated with ${\bm \Upsilon}^1$, ${\bm \Upsilon}^2$,  ${\bm \Upsilon}^3$,  ${\bm \Upsilon}^4$ in the form
\begin{equation*}
{\bf W}(\rho)=
\begin{pmatrix}
Y^1(\rho)    & Y^2(\rho)    & Y^3(\rho)    & Y^4(\rho)    \\[4pt]
(Y^1(\rho))' & (Y^2(\rho))' & (Y^3(\rho))' & (Y^4(\rho))' \\[4pt]
Z^1(\rho)    & Z^2(\rho)    & Z^3(\rho)    & Z^4(\rho)    \\[4pt]
(Z^1(\rho))' & (Z^2(\rho))' & (Z^3(\rho))' & (Z^4(\rho))'
\end{pmatrix} \, ,
\end{equation*}
so that a particular solution $\overline{\bm \Upsilon}=(\overline Y,\overline Z)$ of \eqref{eq:matrix-form} is given by \eqref{eq:Upsilon-bar}.

\subsubsection{The unique solution of \eqref{pbcil3}}
Applying the superposition principle we get \eqref{sol1}.
In order to obtain the unique solution $(Y,Z)$ of the boundary value problem \eqref{pbcil3},  it remains to determine the constants $C_1, C_2, C_3, C_4$ so that the boundary conditions at $\rho= a$ and $\rho= b$ are satisfied.

We check that the constants $C_1, C_2, C_3, C_4$ are uniquely determined. They solve the system
\begin{equation*} \label{eq:sys-A}
A
\begin{pmatrix}
C_1 \\
C_2 \\
C_3 \\
C_4
\end{pmatrix}
=
\begin{pmatrix}
-\left[(\lambda+2\mu)\overline Y'( a)+\tfrac{\lambda}{ a} \, \overline Y( a)+\lambda \frac{\pi k}{h} \,
\overline Z( a)\right] \\[5pt]
-\left[(\lambda+2\mu)\overline Y'( b)+\tfrac{\lambda}{ b} \, \overline Y( b)+\lambda \frac{\pi k}{h} \,
\overline Z( b)\right] \\[5pt]
-\left[\overline Z'( a)-\frac{\pi k}{h} \, \overline Y( a)\right] \\[5pt]
-\left[\overline Z'( b)-\frac{\pi k}{h} \, \overline Y( b)\right]
\end{pmatrix}
\end{equation*}
where the matrix $A=(a_{ij})$, $i,j\in \{1,2,3,4\}$, is given by
\begin{align*}
& a_{1j}=(\lambda+2\mu)(Y^j)'( a)+\tfrac{\lambda}{ a} \, Y^j( a)+\lambda \tfrac{\pi k}{h} \, Z^j( a) \, , \quad a_{2j}=(\lambda+2\mu)(Y^j)'( b)+\tfrac{\lambda}{ b} \, Y^j( b)+\lambda \tfrac{\pi k}{h} \,  Z^j( b) \, ,  \\[5pt]
& a_{3j}=(Z^j)'( a)-\tfrac{\pi k}{h} \, Y^j( a) \, , \quad \hspace{34mm}a_{4j}=(Z^j)'( b)- \tfrac{\pi k}{h} \, Y^j( b) \, .
\end{align*}

We claim that the matrix $A$ is not singular. Consider the homogeneous linear system $A {\bf d}=\0 $ with ${\bf d}=(D_1,D_2,D_3,D_4)^T$. Then
the function ${\bm \Gamma}=(G,H)$ given by
$$
{\bm \Gamma}(\rho)=D_1 {\bm \Upsilon}^1(\rho)+D_2 {\bm \Upsilon}^2(\rho)+D_3 {\bm \Upsilon}^3(\rho)+D_4 {\bm \Upsilon}^4(\rho)
$$
solves system \eqref{eq:hom} coupled with the boundary conditions 
\begin{equation*}
\begin{cases}
(\lambda+2\mu)G'(a)+\frac{\lambda}{a} \, G(a)+\lambda \tfrac{\pi k}{h} H(a)=0  \, ,\\[5pt]
(\lambda+2\mu)G'(b)+\frac{\lambda}{b} \, G(b)+\lambda \tfrac{\pi k}{h} H(b)=0  \, , \\[5pt]
H'(a)-\tfrac{\pi k}{h} G(a)=0  \, , \\[5pt]
H'(b)-\tfrac{\pi k}{h} G(b)=0  \, .
\end{cases}
\end{equation*}
By Lemma \ref{p:prop-ex-un-rad} we then have that ${\bm \Gamma}\equiv (0,0) $ in $(a,b)$ but being ${\bm \Gamma}$ also a solution of system \eqref{eq:hom} for $\rho \in (0,+\infty)$, by local uniqueness for Cauchy problems, ${\bm \Gamma}\equiv (0,0)$ in $(0,+\infty)$. The linear independence of the functions ${\bm \Upsilon}^1$, ${\bm \Upsilon}^2$,
${\bm \Upsilon}^3$, ${\bm \Upsilon}^4$, then implies $D_1=D_2=D_3=D_4=0$.

We just proved that the linear system $A{\bf d}=\0 $ admits only the trivial solution, thus completing the proof.   \endproof

\subsection{Proof of Theorem \ref{teo_soluzioni}} \label{ss:main-teo}

The formal series contained in \eqref{eq:series-expansions} are a consequence of \eqref{eq:comp-expansion}, Lemma \ref{lemma2} and Lemma \ref{p:prop-ex-un-rad}.

It remains to show how those series converge. We start by proving the weak convergence in $H^1(\Omega)$.
Let ${\bm{\mathcal F}}$ be the linear functional defined by ${\bm{\mathcal F}}(\v):=\int_{\partial \Omega} \g \cdot \v \, dS$ for any $\v \in H^1(\Omega;\R^3)$ with $\g$ as in \eqref{eq:def-g}. We observe that thanks to the H\"older inequality and the trace inequality ${\bm{\mathcal F}} \in (H^1(\Omega;\R^3))'$:
\begin{equation*}
   \left|_{(H^1(\Omega;\R^3))'} \langle {\bm{\mathcal F}},\v\rangle_{H^1(\Omega;\R^3)}\right|\le 2p \sqrt{\pi(b^2-a^2)} \,
    C(\Omega) \, \|\v\|_{H^1(\Omega;\R^3)} \qquad \text{for any } \v\in H^1(\Omega;\R^3) \, ,
\end{equation*}
where $C(\Omega)$ is such that $\|{\rm trace}(v)\|_{L^2(\partial\Omega)} \le C(\Omega) \|v\|_{H^1(\Omega)}$ for any $v\in H^1(\Omega)$.

Writing ${\bm{\mathcal F}}=(\mathcal F_1,\mathcal F_2, \mathcal F_3)$ we have that $\mathcal F_1,\mathcal F_2, \mathcal F_3 \in (H^1(\Omega))'$ with $\mathcal F_1=\mathcal F_2$ are the null functionals and
$$
   _{(H^1(\Omega))'} \langle \mathcal F_3,v\rangle_{H^1(\Omega)}=-\int_{C_{a,b}} \chi_p(x,y) \left[v\left(x,y,\tfrac h2\right)-v\left(x,y,-\tfrac h2\right)\right] dxdy \quad \text{for any } v\in H^1(\Omega) \, .
$$
Let us define the sequence of partial sums corresponding to the Fourier expansion in \eqref{eq:distr-expansion}:
\begin{equation*}
  S_M(x,y,z):=\chi_p(x,y) \,  \sum_{m=0}^{M} (-1)^{m+1} \, \frac{4}{h} \sin\left[\frac \pi h (2m+1)z\right] \, .
\end{equation*}

We claim that $S_M \rightharpoonup \mathcal F_3$ weakly in $(H^1(\Omega))'$ as $M\to +\infty$. We first prove that the sequence $\{S_M\}$ is bounded in $(H^1(\Omega))'$. In the next estimate we use the following notations: we put $\widetilde \Omega:=C_{a,b}\times \left(-\frac h2, \frac{3h}{2}\right)$, we still denote by $v$ the symmetric and $2h$-periodic extension of a function $v\in H^1(\Omega)$ (see Section \ref{periodic}) and  by $v\left(x,y,\frac{3h}{2}\right)$ and $v\left(x,y,-\frac h2\right)$, the traces of a function $v\in H^1(\widetilde \Omega)$ on the upper and lower faces of the hollow cylinder $\widetilde \Omega$, respectively.
\begin{align*}
   & \left|_{(H^1(\Omega))'} \langle S_M,v\rangle_{H^1(\Omega)}\right|
   =\left|\sum_{m=0}^{M} (-1)^{m+1} \frac{4}{h} \int_\Omega
     \chi_p(x,y) \sin\left[\frac{\pi}{h}(2m+1)z\right] v(x,y,z) \, dxdydz\right| \\[5pt]
   & \qquad  =\left|\sum_{m=0}^{M} (-1)^{m+1} \frac{2}{h} \int_{\widetilde \Omega} \chi_p(x,y) \sin\left[\frac{\pi}{h}(2m+1)z\right] v(x,y,z) \, dxdydz\right|  \\[7pt]
    & \quad \text{(integration by parts)} \quad =\left|\sum_{m=0}^{M} (-1)^{m+1} \frac{2}{h} \left\{ \int_{C_{a,b}} -\frac{h\cos\left[\frac{3\pi}{2}(2m+1)\right]\chi_p(x,y)}{\pi(2m+1)} \, v\left(x,y,\frac{3h}{2}\right) \, dxdy \right. \right. \\[5pt]
   & \qquad  \qquad \qquad \qquad \qquad \qquad \qquad
     +\int_{C_{a,b}} \frac{h\cos\left[-\frac{\pi}{2}(2m+1)\right]\chi_p(x,y)}{\pi(2m+1)} \, v\left(x,y,-\frac{h}{2}\right) \, dxdy \\[5pt]
   & \qquad \quad \qquad \qquad \qquad \qquad \qquad  \left. \left.
     +\int_{\widetilde \Omega} \frac{h\chi_p(x,y) }{\pi(2m+1)} \, \cos\left[\frac \pi h(2m+1)z\right] \frac{\partial v}{\partial z} (x,y,z)\, dxdydz \right\} \right|
  \end{align*}
  \begin{align*}
   & \quad \text{($2h$-periodicity of $v$)}  \quad =\left|\sum_{m=0}^{M} (-1)^{m+1} \frac{2}{h} \left\{\int_{\widetilde \Omega} \frac{h\chi_p(x,y)}{\pi(2m+1)} \, \cos\left[\frac \pi h(2m+1)z\right] \frac{\partial v}{\partial z} (x,y,z)\, dxdydz \right\} \right| \\[7pt]
   & \qquad =\left|\frac 2 \pi \int_{C_{a,b}} \left\{\sum_{m=0}^{M} \frac{(-1)^{m+1}\chi_p(x,y) }{2m+1}\int_{-\frac h2}^{\frac{3h}{2}}  \cos\left[\frac \pi h(2m+1)z\right] \frac{\partial v}{\partial z} (x,y,z)\,dz \right\}dxdy \right|
  \end{align*}
  \begin{align*}
   & \quad \text{(Cauchy-Schwarz inequality in $\R^{n+1}$)} \\
   & \qquad
    \le \frac{2p} \pi \left[\sum_{m=0}^{M} \frac{1}{(2m+1)^2}\right]^{\frac 12}
    \int_{C_{a,b}} \left[ \sum_{m=0}^{M} \left(\int_{-\frac h2}^{\frac{3h}{2}}  \cos\left[\frac \pi h(2m+1)z\right] \frac{\partial v}{\partial z} (x,y,z)\,dz\right)^2\right]^{\frac 12}  dxdy
 \end{align*}
 \begin{align*}
   & \quad \text{(Bessel inequality)} \quad  \le \frac{2p} \pi \left[\sum_{m=0}^{+\infty} \frac{1}{(2m+1)^2}\right]^{\frac 12}
    \int_{C_{a,b}} \left[ \frac 1h \int_{-\frac h2}^{\frac{3h}{2}} \left(\frac{\partial v}{\partial z} (x,y,z)\right)^2 dz \right]^{\frac 12} dxdy  \\[7pt]
   & \quad \text{(H\"older inequality)} \quad   \le \frac{2p} \pi \left[\sum_{m=0}^{+\infty} \frac{1}{(2m+1)^2}\right]^{\frac 12} \!\!\! \sqrt{\pi(b^2-a^2)}
    \left\{\int_{C_{a,b}} \!\! \left[ \frac 1h \int_{-\frac h2}^{\frac{3h}{2}} \left(\frac{\partial v}{\partial z} (x,y,z)\right)^2 dz \right] dxdy\right\}^{\frac 12}  \\[7pt]
   & \qquad =2p \sqrt{\frac{b^2-a^2}{\pi h}} \left[\sum_{m=0}^{+\infty} \frac{1}{(2m+1)^2}\right]^{\frac 12}
     \left\|\frac{\partial v}{\partial z}\right\|_{L^2(\widetilde \Omega)}\le  4p \sqrt{\frac{b^2-a^2}{\pi h}} \left[\sum_{m=0}^{+\infty} \frac{1}{(2m+1)^2}\right]^{\frac 12}
     \|v\|_{H^1(\Omega)} \, .
\end{align*}
This readily implies
\begin{equation} \label{eq:bound-H1'}
   \|S_M\|_{(H^1(\Omega))'} \le 4 p \sqrt{\frac{b^2-a^2}{\pi h}} \left[\sum_{m=0}^{+\infty} \frac{1}{(2m+1)^2}\right]^{\frac 12}
\end{equation}
and boundedness of $\{S_M\}$ in $(H^1(\Omega))'$ is proved.

Now we claim that
\begin{equation} \label{eq:SM-nuovo}
    \int_\Omega S_M \phi \, d\x \to -\int_{C_{a,b}} \chi_p(x,y) \left[\phi \left(x,y,\tfrac h2\right)-\phi\left(x,y,-\tfrac h2\right)\right] dxdy=_{(H^1(\Omega))'} \langle \mathcal F_3,\phi \rangle_{H^1(\Omega)}
\end{equation}
as $M\to +\infty$, for any $\phi \in C^\infty (\overline \Omega)$. First of all, by using the classical results about pointwise convergence of the Fourier Series applied to suitable $2h$-periodic extensions of the functions
$$
    z\mapsto \phi(x,y,z) \, , \ z\in \left[-\tfrac h 2,0\right] \qquad \text{and} \qquad
     z\mapsto \phi(x,y,z) \, , \ z\in \left[0,\tfrac h 2\right]
$$
one can show that for any $(x,y)\in C_{a,b}$
\begin{equation} \label{eq:pointwise}
   \int_{-\frac h2}^{\frac h2} S_M(x,y,z) \phi(x,y,z) dz \to  -\chi_p(x,y) \left[\phi \left(x,y,\tfrac h2\right)-\phi\left(x,y,-\tfrac h2\right)\right] \, .
\end{equation}
Then applying to the test function $\phi$ the estimates used for proving boundedness of $\{S_M\}$ in $(H^1(\Omega))'$, one can show that for any $M$
\begin{equation}  \label{eq:dominated}
   \left|  \int_{-\frac h2}^{\frac h2} S_M(x,y,z) \phi(x,y,z) dz \right| \le \frac{2\sqrt 2 p}{\pi} \left[\sum_{m=0}^{+\infty} \frac{1}{(2m+1)^2}\right]^{\frac 12} \left\|\frac{\partial \phi}{\partial z}\right\|_{L^\infty(\Omega)} \qquad \text{for any } (x,y)\in C_{a,b} \, .
\end{equation}
By \eqref{eq:pointwise}, \eqref{eq:dominated} and the Dominated Convergence Theorem the proof of \eqref{eq:SM-nuovo} follows. With an essentially similar procedure one can prove that $S_M$ converges in the sense of distributions to $\Lambda_3$ where $\Lambda_3$ is the third component of the vector distribution ${\bm \Lambda}$ defined in \eqref{eq:def-Lambda-1}.

Since $(H^1(\Omega))'$ is a reflexive Banach space, by \eqref{eq:bound-H1'} we infer that along suitable subsequences, the partial sums are weakly convergent in $(H^1(\Omega))'$. Thanks to \eqref{eq:SM-nuovo}, we deduce that the weak limits of this subsequences coincide on the space $C^\infty(\overline\Omega)$ and they equal $\mathcal F_3$ on it. By density of $C^\infty(\overline\Omega)$ in $H^1(\Omega)$, they actually coincide on the whole $H^1(\Omega)$.
This proves that all weakly convergent subsequences weakly converge to $\mathcal F_3$ and hence the sequence $S_M$ is itself weakly convergent to $\mathcal F_3$ in $(H^1(\Omega))'$.
We can now denote by ${\bm S}_M=(0,0,S_M)$ the sequence of vector partial sums in such a way that ${\bm S}_M \rightharpoonup
{\bm{\mathcal F}}$ weakly in $(H^1(\Omega;\R^3))'$ as $M\to +\infty$.

Now, let us consider the linear continuous operator $L$ introduced in the proof of Theorem \ref{p:existence} and its restriction to $V_0^\perp$, where we recall that orthogonality is with respect to the scalar product \eqref{eq:scal-prod-T}.  Then, by Theorem \ref{p:existence} (iv) we deduce that $L_{|V_0^\perp}:V_0^\perp \to (H^1(\Omega;\R^3))'$ is invertible and by the Open Mapping Theorem it follows that its inverse is continuous.

If we define ${\bm U}_M:=L_{|V_0^\perp}^{-1} {\bm S}_M$, then ${\bm U}_M$ is the vector partial sum corresponding to the Fourier expansion \eqref{eq:series-expansions}. Since ${\bm S}_M$ is weakly convergent in $(H^1(\Omega;\R^3))'$ to ${\bm{\mathcal F}}$, then the continuity of $L_{|V_0^\perp}$ implies that ${\bm U}_M$ is weakly convergent in $H^1(\Omega;\R^3)$ to the unique solution $\u$ of \eqref{eq:weak-2} as $M\to +\infty$.

The strong convergence ${\bm U}_M \to \u$ in $L^2(\Omega;\R^3)$ as $M\to +\infty$ is a consequence of the compactness of the embedding $H^1(\Omega;\R^3)\subset L^2(\Omega;\R^3)$.

It remains to prove \eqref{eq:stima-L2}. In order to emphasize the dependence on $k$ we reintroduce it for denoting the functions $Y_k$ and $Z_k$ appearing in the proof of Lemma \ref{p:prop-ex-un-rad}. Testing \eqref{eq:var-rad-0} with $(H,K)=(Y_k,Z_k)$ we have
\begin{equation*}
   \frac{2\mu \pi^2}{h^2} \, k^2 \int_a^b \rho (Z_k(\rho))^2 d\rho \le \int_a^b \rho \Psi_k(\rho) Z_k(\rho) \, d\rho
\end{equation*}
from which we obtain
\begin{equation} \label{eq:L2-Zk}
  \left(\int_a^b (Z_k(\rho))^2 d\rho\right)^{\frac 12}\le \frac{2pb h\sqrt{b-a}}{\mu a \pi^2} \, \frac{1}{k^2} \, .
\end{equation}
Testing again \eqref{eq:var-rad-0} with $(H,K)=(Y_k,Z_k)$ we also have
\begin{equation*}
   2\mu \int_a^b \frac{(Y_k(\rho))^2}{\rho} \, d\rho\le \int_a^b \rho \Psi_k(\rho) Z_k(\rho) \, d\rho
\end{equation*}
from which we obtain
\begin{equation} \label{eq:L2-Yk}
  \int_a^b (Y_k(\rho))^2\, d\rho\le \frac{2pb^2\sqrt{b-a}}{\mu h} \left(\int_a^b (Z_k(\rho))^2 d\rho\right)^{\frac 12}
  \le \frac{4p^2 b^3 (b-a)}{\mu^2 a \pi^2} \frac{1}{k^2}
\end{equation}
where in the last inequality we used \eqref{eq:L2-Zk}.

Let us proceed by considering the difference between the partial sum for $u_1$ and $u_1$ itself:
\begin{align*}
   \|U_M^1-u_1\|_{L^2(\Omega)}^2&=\int_{C_{a,b} } \left\|\sum_{m=M+1}^{+\infty} Y_{2m+1}(\rho) \cos \theta
   \cos\left[\frac{(2m+1)\pi}{h} \, z\right]\right\|_{L^2\left(-\frac h2, \frac h2 \right)}^2 dxdy \\[7pt]
   & =\frac h2 \sum_{m=M+1}^{+\infty} \int_{C_{a,b}} (\cos^2\theta) (Y_{2m+1}(\rho))^2 dxdy
   =\frac{\pi h}{2}\sum_{m=M+1}^{+\infty} \int_a^b (Y_{2m+1}(\rho))^2 \rho \, d\rho \\[7pt]
   &   \le \frac{2 p^2 b^4 (b-a)h}{\mu^2 a\pi} \, \sum_{m=M+1}^{+\infty} \frac{1}{(2m+1)^2}
   \le \frac{p^2 b^4 (b-a)h}{2\mu^2 a\pi} \, \sum_{m=M+1}^{+\infty} \frac{1}{m^2}
   \le \frac{p^2 b^4 (b-a)h}{2\mu^2 a\pi}  \, \frac 1M \, ,
\end{align*}
where we also used \eqref{eq:L2-Yk}.
The estimate for $\|U_n^2-u_2\|_{L^2(\Omega)}$ gives the same result for obvious reasons.

With a completely similar procedure by exploiting this time \eqref{eq:L2-Zk}, we obtain
\begin{align*}
   \|U_n^3-u_3\|_{L^2(\Omega)}^2 \le \frac{2 p^2 b^3 h^3(b-a)}{\mu^2 a^2 \pi^4} \, \sum_{m=M+1}^{+\infty} \frac{1}{(2m+1)^4}
   \le \frac{p^2 b^3 h^3(b-a)}{8\mu^2 a^2 \pi^4} \, \sum_{m=M+1}^{+\infty} \frac{1}{m^4}
   \le \frac{p^2 b^3 h^3(b-a)}{24\mu^2 a^2 \pi^4} \, \frac{1}{M^3} \, .
\end{align*}
The solution we found by means of the Fourier series expansion satisfies \eqref{eq:vanishing-1-2} in the sense of traces of $H^1$-functions.
We conclude the proof of the theorem by observing that this solution coincides with the unique solution of \eqref{eq:weak-2} belonging to $V_0^\perp$. To see this, denote by $\u$ the solution found by means of the Fourier series expansion and by $\w$ the solution in $V_0^\perp$. Both $\u$ and $\w$ possesses the symmetry properties stated in Proposition \ref{prop0} as it occurs to their difference $\u-\w$. But from Theorem \ref{p:existence} we have that $\u-\w\in V_0$ and it is readily seen from \eqref{eq:van-cond} that functions in $V_0$ satisfying those symmetry properties are necessarily the null function.
This proves that $\u=\w$ and completes the proof of the theorem. \endproof

\subsection{Proof of Proposition \ref{p:error}}\label{prooferror}
We rewrite the homogeneous system \eqref{eq:hom} as in \eqref{eq:hom-scaled} so that the corresponding series expansion can be written in the form
\begin{equation}\label{expansion-scaled}
\begin{cases}
{\ds \widetilde Y(t)=\sum_{n=-1}^{+\infty} \widetilde a_n \, t^n+(\ln t) \sum_{n=0}^{+\infty} \widetilde b_n \, t^n} \, , \\[20pt]
{\ds \widetilde Z(t)=\sum_{n=0}^{+\infty} \widetilde c_n \, t^n+(\ln t) \sum_{n=0}^{+\infty} \widetilde d_n \, t^n} \, .
\end{cases}
\end{equation}
The coefficients $\widetilde a_n, \widetilde b_n, \widetilde c_n, \widetilde d_n$ are related to the corresponding coefficients $a_n, b_n, c_n, d_n$ appearing in \eqref{expansion}, by the formulas
\begin{align} \label{eq:coeff-tilde}
& \widetilde a_{-1}=\frac{\pi k}{h} \, a_{-1} \, , \qquad
\widetilde a_n=\left(\frac{h}{\pi k}\right)^n \left[a_n-\ln\left(\frac{\pi k}{h}\right) b_n\right] \, ,\qquad \widetilde b_n=\left(\frac{h}{\pi k}\right)^n  b_n\, , \\[7pt]
& \notag
\widetilde c_n=\left(\frac{h}{\pi k}\right)^n \left[c_n-\ln\left(\frac{\pi k}{h}\right) d_n\right] \, ,\qquad \widetilde d_n=\left(\frac{h}{\pi k}\right)^n d_n\, .
\end{align}
Inserting \eqref{expansion-scaled} into \eqref{eq:hom-scaled} or alternatively combining \eqref{eq:coeff-tilde} and \eqref{eq:linear-system-0}, we see that $\widetilde a_n, \widetilde b_n, \widetilde c_n, \widetilde d_n$ solve the system

\begin{equation} \label{eq:linear-system-scaled}
\begin{cases}
(n^2-1) \widetilde a_n+2n \widetilde b_n+\beta (n-1) \widetilde c_{n-1}+\widetilde \beta \widetilde d_{n-1}=\widetilde\alpha \widetilde a_{n-2} \\
(n^2-1) \widetilde b_n +\widetilde\beta (n-1) \widetilde d_{n-1}=\widetilde\alpha \widetilde b_{n-2}   \\
(n-1)^2 \widetilde c_{n-1}+2(n-1) \widetilde d_{n-1}=\widetilde\delta (n-1) \widetilde a_{n-2}+\widetilde\delta \widetilde b_{n-2}+\widetilde\gamma  \widetilde c_{n-3}  \\
(n-1)^2 \widetilde d_{n-1}=\widetilde\delta  (n-1) \widetilde b_{n-2}+\widetilde\gamma \widetilde d_{n-3}
\end{cases}
\end{equation}
for $n\ge 3$; moreover $\widetilde a_0=\widetilde b_0=\widetilde c_1= \widetilde d_1= \widetilde a_2=\widetilde b_2=0$,
the coefficients $\widetilde a_{-1}, \widetilde a_1, \widetilde b_1, \widetilde c_0$ may be chosen arbitrarily and $\widetilde d_0=\frac{\widetilde\alpha}{\widetilde \beta} \, \widetilde a_{-1}-\frac 2 {\widetilde \beta} \, \widetilde b_1$.

By direct computation one can verify that the unique solution of system \eqref{eq:linear-system-scaled} can be written in form

\begin{equation} \label{eq:Mn}
\begin{pmatrix}
\widetilde a_n \\
\widetilde b_n \\
\widetilde c_{n-1} \\
\widetilde d_{n-1}
\end{pmatrix}
=
\begin{pmatrix}
-\frac{\lambda}{\mu} \, \frac{1}{(n+1)(n-1)} & \frac{2\lambda}{\mu} \, \frac{n}{(n+1)^2(n-1)^2}
& -\frac{\lambda+\mu}{\mu} \, \frac{1}{(n+1)(n-1)^2} & \frac{\lambda+\mu}{\mu} \, \frac{3n+1}{(n+1)^2 (n-1)^3} \\[10pt]
0  &  -\frac{\lambda}{\mu} \, \frac{1}{(n+1)(n-1)}   & 0 & -\frac{\lambda+\mu}{\mu} \, \frac{1}{(n+1)(n-1)^2}  \\[10pt]
\frac{\lambda+\mu}{\mu} \, \frac{1}{n-1}             & -\frac{\lambda+\mu}{\mu} \, \frac{1}{(n-1)^2}
& \frac{\lambda+2\mu}{\mu} \, \frac{1}{(n-1)^2}      & -\frac{2(\lambda+2\mu)}{\mu} \, \frac{1}{(n-1)^3} \\[10pt]
0    &  \frac{\lambda+\mu}{\mu} \, \frac{1}{n-1}     &  0  & \frac{\lambda+2\mu}{\mu} \, \frac{1}{(n-1)^2}
\end{pmatrix}
\begin{pmatrix}
\widetilde a_{n-2} \\
\widetilde b_{n-2} \\
\widetilde c_{n-3} \\
\widetilde d_{n-3}
\end{pmatrix}
\end{equation}

for any $n\ge 3$.

We are interested in the case $n$ odd since when $n$ is even, thanks to \eqref{eq:linear-system-scaled}, we know that $\widetilde a_n=\widetilde b_n=\widetilde c_{n-1}=\widetilde d_{n-1}=0$. Looking at \eqref{eq:Mn}, for any $n\ge 3$ odd, we introduce the matrices
\begin{equation*}
S_n:=\begin{pmatrix}
-\frac{\lambda}{\mu} \, \frac{1}{(n+1)(n-1)} &  -\frac{\lambda+\mu}{\mu} \, \frac{1}{(n+1)(n-1)^2} \\[10pt]
\frac{\lambda+\mu}{\mu} \, \frac{1}{n-1}     & \frac{\lambda+2\mu}{\mu} \, \frac{1}{(n-1)^2}
\end{pmatrix} \, ,
\qquad
T_n:=\begin{pmatrix}
\frac{2\lambda}{\mu} \, \frac{n}{(n+1)^2 (n-1)^2} &  \frac{\lambda+\mu}{\mu} \, \frac{3n+1}{(n+1)^2 (n-1)^3} \\[10pt]
-\frac{\lambda+\mu}{\mu} \, \frac{1}{(n-1)^2}     & -\frac{2(\lambda+2\mu)}{\mu} \, \frac{1}{(n-1)^3}
\end{pmatrix} \, .
\end{equation*}
In this way, system \eqref{eq:Mn} may written in the form
\begin{align*} 
&
\begin{pmatrix}
\widetilde b_n \\[2pt]
\widetilde d_{n-1}
\end{pmatrix}
=S_n
\begin{pmatrix}
\widetilde b_{n-2} \\[2pt]
\widetilde d_{n-3}
\end{pmatrix} \, , \qquad
\begin{pmatrix}
\widetilde a_n \\[2pt]
\widetilde c_{n-1}
\end{pmatrix}
=S_n
\begin{pmatrix}
\widetilde a_{n-2} \\[2pt]
\widetilde c_{n-3}
\end{pmatrix}
-T_n
\begin{pmatrix}
\widetilde b_{n-2} \\[2pt]
\widetilde d_{n-3}
\end{pmatrix} \, .
\end{align*}
After an iterative procedure we may write
\begin{align} \label{eq:iterative-scheme}
\begin{cases}
\begin{pmatrix}
\widetilde b_n \\[2pt]
\widetilde d_{n-1}
\end{pmatrix}
=\left(\ds{\prod_{m=0}^{(n-3)/2}} S_{n-2m}\right)
\begin{pmatrix}
\widetilde b_1 \\[2pt]
\widetilde d_{0}
\end{pmatrix} \, , \\[15pt]
\begin{pmatrix}
\widetilde a_n \\[2pt]
\widetilde c_{n-1}
\end{pmatrix}
=\left(\ds{\prod_{m=0}^{(n-3)/2}} S_{n-2m}\right)
\begin{pmatrix}
\widetilde a_1 \\[2pt]
\widetilde c_{0}
\end{pmatrix}
-\ds{\sum_{j=0}^{(n-3)/2}} \left[ \left(\ds{\prod_{m=1}^{j} S_{n-2m+2}}\right) T_{n-2j} \left(\ds{\prod_{m=j+1}^{(n-3)/2} S_{n-2m}}\right)
\begin{pmatrix}
\widetilde b_1 \\[2pt]
\widetilde d_{0}
\end{pmatrix} \right]  ,
\end{cases}
\end{align}
for any $n\ge 3$ odd, with the convention that for any sequence of matrices $A_m\in \R^{2\times 2}$
\begin{equation*}
\prod_{m=m_1}^{m_2} A_m=
\begin{pmatrix}
1 & 0 \\
0 & 1
\end{pmatrix} \qquad \text{and} \qquad \sum_{m=m_1}^{m_2} A_m
=
\begin{pmatrix}
0 & 0 \\
0 & 0
\end{pmatrix}
\end{equation*}
whenever $m_1>m_2$.

By induction one can verify that for any $j\le \frac{n-3}{2}$
\begin{align*}
& \prod_{m=0}^{j} S_{n-2m}=
\begin{pmatrix}
-\frac{(j+1)\lambda+j \mu}{\mu (n+1) [n+1-2(j+1)] \prod_{m=1}^{j} (n+1-2m)^2 } &   \qquad -\frac{(j+1)(\lambda+\mu)}{\mu (n+1) \prod_{m=1}^{j+1} (n+1-2m)^2 }   \\[15pt]
\frac{(j+1)(\lambda+\mu)}{\mu [n+1-2(j+1)] \prod_{m=1}^{j} (n+1-2m)^2 }   &  \qquad
\frac{(j+1)\lambda+(j+2)\mu}{\mu \prod_{m=1}^{j+1} (n+1-2m)^2}
\end{pmatrix}
\end{align*}

and, in turn, by \eqref{eq:characterization-infty} we infer

\begin{equation} \label{eq:prod-Sn}
\left\|\prod_{m=0}^{j} S_{n-2m}\right\|_\infty \le \frac{(\lambda+\mu) (n-2j)(j+2)}
{\mu \prod_{m=1}^{j+1} (n+1-2m)^2 } \, .
\end{equation}

In particular, with appropriate choices of the minimum and the maximum values of the index in the product \eqref{eq:prod-Sn} and with appropriate changes of index, for any $n\ge 3$ odd, we obtain the estimates

\begin{align} \label{eq:prod-Sn-bis}
& \left\|\prod_{m=0}^{(n-3)/2} S_{n-2m}\right\|_\infty \le \frac{3(\lambda+\mu) (n+1)}{\mu2^{n} \left[\left(\frac{n-1}2\right)!\right]^2} \, , \qquad
\left\|\prod_{m=1}^{j} S_{n-2m+2}\right\|_\infty \le \frac{(\lambda+\mu) (n-2j+2)(j+1)}
{\mu \prod_{m=1}^{j} (n+1-2m)^2 } \, , \\[15pt]
\notag & \left\|\prod_{m=j+1}^{(n-3)/2} S_{n-2m}\right\|_\infty \le
\frac{3(\lambda+\mu)  (n-2j-1)}{2\mu \prod_{m=j+2}^{(n-1)/2} \, (n+1-2m)^2} \, .
\end{align}

On the other hand, we observe that for the components of the matrices $S_n$ and $T_n$ the following inequalities hold true:
\begin{equation*} 
|(T_n)_{ij}|\le \tfrac{3}{n-1} \, |(S_n)_{ij}| \qquad \text{for any $i,j \in \{1,2\}$ and $n\ge 3$,}
\end{equation*}
which, in turn, implies $\|T_n\|_\infty \le \frac{3}{n-1} \, \|S_n\|_\infty=3(\lambda+\mu)\frac{ 2n}
{\mu (n-1)^3 }$; the last inequality is obtained by \eqref{eq:prod-Sn} with $j=0$.

Therefore, combining \eqref{eq:prod-Sn} and \eqref{eq:prod-Sn-bis}, for any $n\ge 3$ odd, we obtain
\begin{align} \label{eq:prod-Sn-ter}
& \left\|\ds{\sum_{j=0}^{(n-3)/2}} \left[ \left(\ds{\prod_{m=1}^{j} S_{n-2m+2}}\right) T_{n-2j} \left(\ds{\prod_{m=j+1}^{(n-3)/2} S_{n-2m}}\right)
\right] \right\|_\infty  \\[5pt]
\notag  & \ \ \ \le \ds{\sum_{j=0}^{(n-3)/2}} \left\|\ds{\prod_{m=1}^{j} S_{n-2m+2}}\right\|_\infty  \left\|T_{n-2j} \right\|_\infty \left\|\ds{\prod_{m=j+1}^{(n-3)/2} S_{n-2m}}\right\|_\infty \\[5pt]
\notag &    \le  \ds{\sum_{j=0}^{(n-3)/2}} \frac{18(\lambda+\mu)^3 (n-2j+2)(n-2j)(j+1)}{\mu^3 \,  2^{n} \left[\left(\frac{n-1}{2}\right)!\right]^2}
\le \frac{9 (\lambda+\mu)^3 \, n(n+2)(n^2-1)}{4\mu^3 \,  2^{n} \left[\left(\frac{n-1}{2}\right)!\right]^2} \, ,
\end{align}
where in the last inequality we used the estimate $(n-2j+2)(n-2j)\le n(n+2)$ and the identity $\sum_{j=0}^{(n-3)/2} (j+1)=\frac{n^2-1}{8}$.

Combining \eqref{eq:norm-infty} with  \eqref{eq:iterative-scheme}, \eqref{eq:prod-Sn-bis} and \eqref{eq:prod-Sn-ter}, for any $n\ge 3$ odd, we obtain
\begin{align} \label{eq:final-1-bis}
\left| \begin{pmatrix}
\widetilde a_n \\
\widetilde c_{n-1}
\end{pmatrix}
\right|_\infty & \le \frac{3(\lambda+\mu) (n+1)}{\mu 2^n \left[\left(\frac{n-1}2\right)!\right]^2}
\left|
\begin{pmatrix}
\widetilde a_1 \\
\widetilde c_0
\end{pmatrix}
\right|_\infty
+\frac{9(\lambda+\mu)^3 \, n(n+2)(n^2-1)}{4 \mu^3 \,  2^{n} \left[\left(\frac{n-1}{2}\right)!\right]^2}
\left|
\begin{pmatrix}
\widetilde b_1 \\
\widetilde d_0
\end{pmatrix}
\right|_\infty \\[10pt]
\notag &  \le \frac{3(2\lambda+5\mu)(\lambda+\mu)^2 \, (n+1)(3n^3+3n^2-6n+4)}{4 \mu^3 \,  2^{n} \left[\left(\frac{n-1}{2}\right)!\right]^2} \, \max\{\widetilde a_{-1}, \widetilde a_1, \widetilde b_1, \widetilde c_0\}
\end{align}
and
\begin{equation} \label{eq:final-2-bis}
\left| \begin{pmatrix}
\widetilde b_n \\
\widetilde d_{n-1}
\end{pmatrix}
\right|_\infty \le  \frac{3(\lambda+\mu) (n+1)}{\mu 2^n \left[\left(\frac{n-1}2\right)!\right]^2}
\left|
\begin{pmatrix}
\widetilde b_1 \\
\widetilde d_0
\end{pmatrix}
\right|_\infty
\le \frac{3(2\lambda+5\mu) (n+1)}{\mu 2^n \left[\left(\frac{n-1}2\right)!\right]^2} \, \max\{\widetilde a_{-1}, \widetilde a_1, \widetilde b_1, \widetilde c_0\}
\end{equation}
where we exploited the fact that $\widetilde d_0=\frac{\widetilde\alpha}{\widetilde \beta} \, \widetilde a_{-1}-\frac 2 {\widetilde \beta} \, \widetilde b_1$, accordingly with what already explained in the lines below \eqref{eq:linear-system-scaled}, so that
\begin{equation*}
|\widetilde d_0|\le \frac{\widetilde\alpha+2}{\widetilde\beta} \max\{\widetilde a_{-1}, \widetilde b_1\}= \frac{2\lambda+5\mu}{\lambda+\mu} \max\{\widetilde a_{-1}, \widetilde b_1\} \, ,
\end{equation*}
from which it follows that
\begin{equation*}
\max\left\{\left| \begin{pmatrix}
\widetilde a_1 \\
\widetilde c_0
\end{pmatrix}
\right|_\infty ,
\left| \begin{pmatrix}
\widetilde b_1 \\
\widetilde d_0
\end{pmatrix}
\right|_\infty \right\} \le \frac{2\lambda+5\mu}{\lambda+\mu}
\max\{\widetilde a_{-1}, \widetilde a_1, \widetilde b_1, \widetilde c_0\}  \, .
\end{equation*}

Since we are interested to the restrictions of the functions $Y$ and $Z$ to the interval $[a,b]$, we have to evaluate the series expansion \eqref{expansion-scaled} of the functions $\widetilde Y$ and $\widetilde Z$ for
$t\in \left[\frac{\pi k}{h} a, \frac{\pi k}{h} b\right]$.

Let $N$ odd be the number at which we want to truncate the series expansions in \eqref{expansion-scaled}. Recalling that the coefficients $\widetilde a_n, \widetilde b_n, \widetilde c_{n-1}, \widetilde d_{n-1}$ vanish for $n$ even, we may write
\begin{equation*} 
\begin{cases}
{\ds \widetilde Y(t)=\left(\, \sum_{n=-1}^{N} \widetilde a_n \, t^n+(\ln t) \sum_{n=0}^{N} \widetilde b_n \, t^n\right) + \left(\, \sum_{n=N+2}^{+\infty} \widetilde a_n \, t^n+(\ln t) \sum_{n=N+2}^{+\infty} \widetilde b_n \, t^n\right) } \, , \\[20pt]
{\ds \widetilde Z(t)=\left(\, \sum_{n=0}^{N-1} \widetilde c_n \, t^n+(\ln t) \sum_{n=0}^{N-1} \widetilde d_n \, t^n \right)
	+\left(\, \sum_{n=N+1}^{+\infty} \widetilde c_n \, t^n+(\ln t) \sum_{n=N+1}^{+\infty} \widetilde d_n \, t^n\right)} \, ,
\end{cases}
\end{equation*}
and define the truncation error as
\begin{equation*}
E_{k,N}=\max \left\{\max\limits_{t\in \left[\frac{\pi k}{h} a, \frac{\pi k}{h} b\right]} \left|\, \sum_{n=N+2}^{+\infty} \widetilde a_n \, t^n+(\ln t) \sum_{n=N+2}^{+\infty} \widetilde b_n \, t^n\right| , \,
\max\limits_{t\in \left[\frac{\pi k}{h} a, \frac{\pi k}{h} b\right]} \left|\, \sum_{n=N+1}^{+\infty} \widetilde c_n \, t^n+(\ln t) \sum_{n=N+1}^{+\infty} \widetilde d_n \, t^n\right|  \right\}
\end{equation*}
By \eqref{eq:final-1-bis} and \eqref{eq:final-2-bis}, we see that for any $t\in \left[\frac{\pi k}{h} a, \frac{\pi k}{h} b\right]$ we have
\begin{align*}
0& \le E_{k,N}\le \widetilde C(a,b,k) \,
\left[ \, \sum_{n=N+2}^{+\infty} \left(\frac{\pi k b}{h}\right)^n
\left|
\begin{pmatrix}
\widetilde a_n \\[5pt]
\widetilde c_{n-1}
\end{pmatrix}
\right|_\infty + \sum_{n=N+2}^{+\infty} \left(\frac{\pi k b }{h}\right)^n
\left|
\begin{pmatrix}
\widetilde b_n \\[5pt]
\widetilde d_{n-1}
\end{pmatrix}
\right|_\infty \, \right]
\\[5pt]
&  \le  \widetilde C(a,b,k) \max\{\widetilde a_{-1}, \widetilde a_1, \widetilde b_1, \widetilde c_0\} \sum_{\underset{\text{$n$ odd}}{n=N+2}}^{+\infty} \left(\frac{\pi k b }{h}\right)^n \frac{3(2\lambda+5\mu)(\lambda+\mu)^2 \,(n+1) (3n^3+3n^2-6n+8)}{4\mu^3 \,  2^{n} \left[\left(\frac{n-1}{2}\right)!\right]^2}
\end{align*}
where we put $\widetilde C(a,b,k)=\max\left\{1,\frac{h}{\pi kb}\right\}  \max\Big\{1,|\ln\left(\frac{\pi k a }{h}\right)|,|\ln\left(\frac{\pi k b }{h}\right)|\Big\}$.

Since we are interested to truncation of the series expansion with a sufficiently large number of terms, letting $P(n):=(n+1) (3n^3+3n^2-6n+8)$, it is not restrictive to assume $N \ge 3$ in such a way that the sequence $n\mapsto 2^{-n}P(n)$ becomes decreasing for $n\ge N+2\ge 5$.

In this way, for $N \ge 3$ odd, we obtain for all $t\in \left[\frac{\pi k}{h} a, \frac{\pi k}{h} b\right]$
\begin{align} \label{eq:E-n(t)}
0& \le E_{k,N}\le \widetilde C(a,b,k) \max\{\widetilde a_{-1}, \widetilde a_1, \widetilde b_1, \widetilde c_0\} \frac{3(2\lambda+5\mu)(\lambda+\mu)^2 \, P(N+2)}{4\mu^3 \,  2^{N+2} \left(\frac{N+1}{2}\right)!} \sum_{m=\frac{N+1}{2}}^{+\infty} \frac{\left(\frac{\pi k b}{h}\right)^{2m+1}}{m!}  \\[10pt]
\notag  & \le \widetilde C(a,b,k) \max\{\widetilde a_{-1}, \widetilde a_1, \widetilde b_1, \widetilde c_0\}
\left(\frac{\pi k b}{h}\right)^{N+2} e^{\left(\frac{\pi k b}{h}\right)^2} \,
\frac{3(2\lambda+5\mu)(\lambda+\mu)^2 \, P(N+2)}{16\mu^3 \,  2^{N}
 \left[\left(\frac{ N+1}{2}\right)!\right]^2}  \, ,
\end{align}
where in the last estimate we used the Lagrange form of the reminder in the Taylor formula for the exponential function and
$P(N+2)=(N+3)(3N^3+21N^2+42N +32)$.

According to the rescaling introduced in \eqref{eq:hom-scaled} one may define the functions $\widetilde{\bm \Upsilon}^j$, whose series expansions are given by \eqref{expansion-scaled}
with coefficients in \eqref{eq:coeff-tilde} and with $a_{-1}, a_1, b_1, c_0$ given by \eqref{eq:cases} in the cases corresponding to $j\in \{1,2,3,4\}$.

In these four cases, the quantity $\max\{\widetilde a_{-1}, \widetilde a_1, \widetilde b_1, \widetilde c_0\}$, appearing in the right hand side of \eqref{eq:E-n(t)}, admits the following estimates:
\begin{equation} \label{eq:est-max}
\begin{cases}
\max\{\widetilde a_{-1}, \widetilde a_1, \widetilde b_1, \widetilde c_0\}= \frac{\pi k}h \, \max\left\{1,\frac{\mu}{\lambda+\mu} \, \ln\left(\frac{\pi k}h\right) \right\} \,  & \qquad \text{if $j=1$,} \\[7pt]
\max\{\widetilde a_{-1}, \widetilde a_1, \widetilde b_1, \widetilde c_0\}=\frac h {\pi k} & \qquad \text{if $j=2$,} \\[7pt]
\max\{\widetilde a_{-1}, \widetilde a_1, \widetilde b_1, \widetilde c_0\}=\frac h {\pi k} \,
\max\left\{1,\frac{2(\lambda+2\mu)}{\lambda+\mu} \, \ln\left(\frac{\pi k}h\right) \right\} & \qquad \text{if $j=3$,} \\[7pt]
\max\{\widetilde a_{-1}, \widetilde a_1, \widetilde b_1, \widetilde c_0\}=1 & \qquad \text{if $j=4$.} \\
\end{cases}
\end{equation}
For $k>1$ is easy to see that all the maximum in \eqref{eq:est-max} are less or equal than
$$
  \max\Big\{\tfrac{\pi k}h ,\tfrac{\mu}{\lambda+\mu}\, \tfrac{\pi k}h \ln\left(\tfrac{\pi k}h \right), \tfrac{2(\lambda+2\mu)}{\lambda+\mu} \, \tfrac h {\pi k} \ln\left(\tfrac{\pi k}h \right) \Big\} \, ,
$$
so that \eqref{err} follows. \endproof

\section{Conclusions}\label{conclusions}
In this work we started from an applicative problem, suggested by Studio De Miranda Associati, an engineering company expertized in building long span bridges. They proposed to study the blister, a structural element in bridges where the steel forestay anchors to the deck. The aim is to obtain an explicit formula to estimate the tensions in the blister, useful for the practical design of bridges.

The problem can be solved through the resolution of the elasticity equation with a specific geometry and load configuration. Hence, the first step was to define the geometry of the element. Through some simplifications we end up with a hollow circular cylinder axially loaded at the end faces; the volume of the cylinder represents the portion of the deck concrete where the stresses diffusion happens, while the applied load is given by the force that the stay has to transfer to the deck.   Clearly this geometry and load configuration can be refined in order to model a real blister, but this is a first step in this way and we leave more sophisticated models to future works.

As matter of fact, from literature we learn that the elasticity equation was explicitly solved only for very particular domains and load conditions, e.g. in prisms \cite{sundara}. 
In this paper we provide the explicit solution for the hollow cylinder axially loaded, proceeding by steps: first of all we provide a periodic extension of the load in $z$ direction, so that we expand the solution in Fourier series with respect to the variable $z$. Then we compute the Fourier coefficients in $x$ and $y$ passing to cylindrical coordinates and expanding such functions in power series. In Theorem \ref{teo_soluzioni} we write the explicit solution for the problem, written in series expansion. We point out that this solution may have an own interest in the construction science field, beyond the application to the blister. 

To employ directly the formula in real situations, such as the blister design, it is necessary to consider approximated solutions, giving some estimates on the errors due to the truncating of the series. In Section \ref{num} we proposed a case of study, where, fixing the parameters involved in the problem, we are able to find the distribution of the stresses in the cylinder. These plots can be obtained through a simple code, written in MATLAB$^\circledR$ or GNU Octave$^\circledR$, running in brief time, e.g. 1-3 minutes, depending on the number where we truncate the series.

From these results it is possible to find the maximum and the minimum of the different stresses acting on the cylinder, their position on the element and an estimate on the error due to the truncation of the series. Knowing these values, the engineering designer can choice for instance the most appropriate strand anchorage from the commercial catalogue, see Figure \ref{fig:detail}, in order to not exceed specific limit stresses in the reinforced concrete. Since the map of the tensions is given, see e.g. Figure \ref{plots}, the engineer can design the steel reinforcements in the concrete, at least on a pre-dimensioning level, and can check the concrete cracking stresses. 

As we explained, to get more precise results on realistic blisters we should modify the geometry of the element and  the configuration of the loads; this may be a future work, but we point out that, more the geometry and the distribution of the loads are complex more the expectations to find explicit solutions are few, so that the finite element analysis may be preferred.
\bigskip

{\bf Notations} 
We give some notations that will be used throughout this paper about functional spaces and differential operators acting on scalar functions, vector valued functions, matrix valued functions. We denote by $\Omega$ a general domain in $\R^N$, $N\ge 1$ where by domain we mean a connected open set in $\R^N$.

\begin{itemize}
	
	\item Given two vectors $\x=(x_1, \dots, x_N), \y=(y_1,\dots, y_N) \in \R^N$ we denote by $\x \cdot \y=\sum_{i=1}^{N} x_i y_i$ their Euclidean scalar product and by $|\x|=\sqrt{\x\cdot \x}$ the Euclidean modulus of $\x$;
	
	\item the $\infty$-norm of vectors is	$|\x|_\infty :=\max\limits_{1\le i \le N} \, |x_i|$;
	
	\item $\R^{M \times N}$: space of $M \times N$ matrices;

	\item if $A\in \R^{M \times N}$ and $\x \in \R^N$ is a vector, $A\x$ denotes the usual product of matrices where $\x$ has to be seen as a vector column;
	
	\item letting $A=(a_{ij}), B=(b_{ij}) \in \R^{N \times N}$ we denote by $A:B=\sum_{i,j=1}^{N} a_{ij} b_{ij}$ their Euclidean scalar product and by $|A|=\sqrt{A:A}$ its Euclidean modulus;
	
	\item given $A\in \R^{M \times N}$ we denote by $A^T\in \R^{N \times M}$ its transpose;
	
	\item given $A\in \R^{N \times N}$ we introduce the operator $\infty$-norm of matrices by
	$\|A\|_\infty:=\sup\limits_{\x\in \R^N \setminus \{ \0 \} } \, \frac{|A\x|_\infty}{|x|_\infty}$
	so that we have in particular
	\begin{equation} \label{eq:norm-infty}
	|A\x|_\infty \le \|A\|_\infty \, |\x|_\infty  \qquad \text{for any } \x\in \R^N.
	\end{equation}
Letting $A=(a_{ij}),\in \R^{N}$, the following characterization of $\|\cdot \|_\infty$ holds:
	\begin{equation}  \label{eq:characterization-infty}
	\|A\|_\infty=\max_{1\le i\le N} \, \sum_{j=1}^{N} |a_{ij}|;
	\end{equation}
	being $\|\cdot \|_\infty$ an operator norm, it is \textit{sub-multiplicative} in the sense that
$\|AB\|_\infty \le \|A\|_\infty \, \|B\|_\infty$ for any  $A,B\in \R^{N \times N}$.
	
	\item some well known functional spaces of functions defined from on an open set $\Omega\subset \R^N$ to a vector space $V$ which could be $\R^M$ or a space of matrices: $C^k(\Omega;V)$, $L^p(\Omega;V)$, $H^k(\Omega;V)$ with $0\le k\le \infty$ integer and $1\le p \le \infty$;
	
	\item for $0\le k\le \infty$ integer, $C^k(\overline\Omega;V)$ denotes the space of restrictions to $\overline \Omega$ of functions in $C^k(\R^N;V)$;
	
	\item $\mathcal D(\Omega;V)$: space of $C^\infty(\Omega;V)$ with compact support in $\Omega$;
	
	\item $\mathcal D'(\Omega;V)$: space of vector distributions, i.e. the dual space of $\mathcal D(\Omega;V)$;
	
	\item given a scalar function $g\in C^1(\Omega;\R)$, we denote by $\nabla g\in C^0(\Omega;\R^n)$ its gradient;
	
	\item given a vector valued function $\u\in C^1(\Omega;\R^M)$, we denote by $\nabla \u\in C^0(\Omega;\R^{M \times N})$ its Jacobian matrix;
	
	\item given a vector valued function $\u\in C^1(\Omega;\R^N)$, $\Omega \subseteq \R^N$, we denote by $\D \u\in C^0(\Omega;\R^{N \times N})$ its symmetric gradient defined by $\D \u=\dfrac{\nabla \u+\nabla \u^T}{2}$ (linearized strain tensor when $N=3$);
	
	\medskip
	
	\item given $U \in C^1(\Omega; \R^{M \times N})$, $\Omega \subseteq \R^N$, we denote by $\dive \, U \in C^0(\Omega;\R^M)$ the vector field $\v=(v_1,\dots,v_M)$ such that $v_i=\sum_{j=1}^{N} \frac{\partial U_{ij}}{\partial x_j}$, $i=1,\dots, M$;
	
	\item given $\u=(u_1,\dots,u_M) \in C^2(\Omega; \R^M)$, we denote by $\Delta \u\in C^0(\Omega; \R^M)$ the Laplacian of $\u$ defined component by component, i.e. $\Delta \u=(\Delta u_1, \dots, \Delta u_M)$ where in the last identity $\Delta$ denotes the usual Laplacian of a real valued function.
\end{itemize}

\bigskip

{\bf Acknowledgments} The two authors are members of the Gruppo Nazionale per l'Analisi Matematica, la Probabilit\`{a} e le loro Applicazioni (GNAMPA) of the Istituto Nazionale di Alta Matematica (INdAM). The second author acknowledges partial financial support from the PRIN project 2017 ``Direct and inverse problems for partial differential equations: theoretical aspects and applications''. The authors acknowledge partial financial support from the INdAM - GNAMPA project 2022 ``Modelli del 4° ordine per la dinamica di strutture ingegneristiche: aspetti analitici e applicazioni''.

The second author acknowledges partial financial support from the research project ``Metodi e modelli per la matematica e le sue applicazioni alle scienze, alla tecnologia e alla formazione'' Progetto di Ateneo 2019 of the University of Piemonte Orientale ``Amedeo Avogadro''.

\end{document}